\documentclass[11pt]{article}

\usepackage{array}
\usepackage{graphicx}
\usepackage{amsmath}
\usepackage{amssymb}
\usepackage{hhline}
\usepackage{bm}
\usepackage{cancel}

\title{\bf Stability analysis of the \\
numerical Method of characteristics \\
applied to a class of \\
energy-preserving hyperbolic systems. \\
Part I: Periodic boundary conditions}

\author{ T.I. Lakoba\footnote{tlakoba@uvm.edu, \ 1 (802) 656-2610}, \ Z. Deng
 \vspace{0.5cm} \\
  Department of Mathematics and Statistics, 16 Colchester Ave., \\
 University of Vermont, Burlington, VT 05401, USA}

\newcommand{\noi}{\noindent}
\newcommand{\und}{\underline}
\newcommand{\D}{\Delta}
\newcommand{\be}{\begin{equation}}
\newcommand{\ee}{\end{equation}}
\newcommand{\bsube}{\begin{subequations}}
\newcommand{\esube}{\end{subequations}}
\newcommand{\ba}{\begin{array}}
\newcommand{\ea}{\end{array}}
\newcommand{\To}{\rightarrow}

\newcommand{\bea}{\begin{eqnarray}}
\newcommand{\eea}{\end{eqnarray}}
\newcommand{\so}{\Rightarrow}

\newcommand{\sech}{\mbox{sech}}

\newcommand{\dt}{\D t}

\newcommand{\Sp}{S^{+}}
\newcommand{\Sm}{S^{-}}
\newcommand{\Spm}{S^{\pm}}

\newcommand{\vSp}{{\bf \und{S}}^{+}}
\newcommand{\vSm}{{\bf \und{S}}^{-}}
\newcommand{\vSpm}{{\bf \und{S}}^{\pm}}
\newcommand{\vSmp}{{\bf \und{S}}^{\mp}}
\newcommand{\tSp}{s^{+}}
\newcommand{\tSm}{s^{-}}
\newcommand{\tSpm}{s^{\pm}}
\newcommand{\tvS}{{\bf s}}
\newcommand{\tvSp}{\und{s}^{+}}
\newcommand{\tvSm}{\und{s}^{-}}
\newcommand{\tvSpm}{\und{s}^{\pm}}

\begin{document}
\baselineskip 18 pt

\maketitle

\vspace*{2cm}

\begin{center}
 {\bf Abstract}
\end{center}

We study numerical (in)stability of the Method of characteristics (MoC)
applied to a system of non-dissipative hyperbolic partial differential equations
(PDEs) with periodic boundary conditions. We consider three different solvers
along the characteristics: simple Euler (SE), modified Euler (ME), and
Leap-frog (LF). The two former solvers are well known to exhibit a mild, but
unconditional, numerical instability for non-dissipative ordinary differential
equations (ODEs). They are found to have a similar (or stronger,
for the MoC-ME) instability when applied to non-dissipative PDEs. 
On the other hand, the LF
solver is known to be stable when applied to non-dissipative ODEs.
However, when applied to non-dissipative PDEs within the MoC framework,
it was found to have by far the strongest instability among all three solvers.
We also comment on the use of the fourth-order Runge--Kutta solver within the
MoC framework.

\vskip 1.1 cm

\noi

\bigskip


\newpage

\section{Introduction}

In this series of two papers we address numerical stability of the Method
of characteristics (MoC) applied to a class of non-dissipative hyperbolic
partial differential equations (PDEs) in 1 time + 1 space dimension. For
reference purposes, we will now present the idea of the MoC using the
following system as an example:
\be
\ba{rcl}
u_{1t} + u_{1x} & = & f_1(u_1,u_2)  \\
u_{2t} - u_{2x} & = & f_2(u_1,u_2),  
\ea
\label{e01}
\ee
where $f_{1,2}$ are some functions. In the MoC, each of the equations in
\eqref{e01} is transformed to an ordinary differential equation (ODE) in its
respective variable, $\eta_{\pm}=x\mp t$, and then is solved by an ODE
numerical solver.

The MoC is widely used and is described in most textbooks on numerical solution
of PDEs. It is, therefore, quite surprising that its stability has not been
investigated in as much detail as that of most other numerical methods.
For example, in \cite{GKO}, it was considered by the von Neumann
analysis (which implies periodic boundary conditions (BC)) for a simple
scalar model
\bsube
\be
u_t+u_x = 0,
\label{e02a}
\ee
with the ODE solver being the simple Euler method. Clearly, that analysis
could be straightforwardly generalized for a vector model
\be
{\bf u}_t+{\bf A u}_x = {\bf 0},
\label{e02b}
\ee
\label{e02}
\esube
where ${\bf A}$ is any constant diagonalizable matrix with real eigenvalues; 
an analysis for a
similar model in three spatial dimensions was presented in \cite{stab_3DMoC}.
However, we have been unable to find a systematic stability analysis ---
even for periodic BC --- for a system of the form
\be
{\bf u}_t+{\bf A u}_x = {\bf B u},
\label{e03}
\ee
where ${\bf A,\,B}$ are constant matrices. A partial exception is a conference
paper \cite{Australian08}, where a model with a somewhat more complicated
right-hand side (rhs), describing a specific engineering application, was
studied by the von Neumann analysis. The focus of that paper was on the
examination of the effect of various parameters of that particular model rather
than on the impact on stability of the MoC by the ODE solver used. (In fact,
the ODE solver used in \cite{Australian08} was the simple Euler method.)

In this paper we will examine the effect of the ODE solver on the stability of
the MoC applied to the model problem \eqref{e03} with {\em periodic} BC.
We will further limit the scope of the problem as follows. First, and most
importantly, we will consider only non-dissipative and stable systems. This
implies that matrix ${\bf B}$ possesses only imaginary eigenvalues. Second,
we will consider only three ODE solvers: simple Euler (SE), modified Euler (ME),
and the Leap-frog (LF). Let us note that for a non-dissipative stable ODE,
the numerical solutions obtained by the SE and ME are known to be mildly,
but unconditionally, unstable, with the growth rates of the numerical instability
being $O(\dt)$ and $O(\dt^3)$, respectively; see, e.g., Secs.~3 and 4 below. 
On the contrary, the LF solver is known to
quasi-preserve (i.e., not shift systematically) at least some of the conserved
quantities of non-dissipative ODEs.

Our third assumption is mostly cosmetic and does not affect the methodology 
of our analysis. Namely, we will assume that the constant matrix ${\bf A}$
in \eqref{e03} can be diagonalized into the form
\bsube
\be
{\bf A} = c_1{\bf I}_N + c_2{\bf \Sigma}, \qquad 
{\bf \Sigma} = {\rm diag}({\bf I}_{N_1},-{\bf I}_{N_2}),
\label{e04a}
\ee
where $c_{1,2}$ are scalars, ${\bf I}_N$ is the $N\times N$ identity matrix, and
$N$ is the dimension of vector ${\bf u}$. Without loss of generality (i.e.,
by a simple change of variables) one can set 
\be
c_1=0, \quad c_2=1.
\label{e04b}
\ee
\label{e04}
\esube
We will also let $N_1=N_2=N/2$. Particular conclusions of the forthcoming analysis
may change if the diagonalization of matrix ${\bf A}$ is different from \eqref{e04a}.
However, the methodology of the analysis will not be affected.

Our analysis reveals two facts which, to our knowledge, have not been previously
pointed out. First, in contrast to the situation with most numerical methods,
the most numerically unstable Fourier harmonics may occur not at the edges of the
spectrum but in the ``middle" (see the footnote in Section 4) 
of it. In fact, the (in)stability of the highest
Fourier harmonics is the same as that of the lowest ones. This fact easily follows
from the foregoing analysis and holds for the MoC employing any ODE solver
as long as matrix ${\bf A}$ has the form \eqref{e04}.

Second, and quite unexpectedly, 
the LF, which outperforms the SE and ME when applied to non-dissipative
ODEs, performs much worse than those methods when applied to non-dissipative
stable PDEs within the MoC framework.

The main part of this work is organized as follows. In Section 2
we specify the physical model to be considered below.
Our analysis will be developed for the linearized version of that model.
In Section 3 we present the von Neumann
analysis and its verification by direct numerical simulations of the PDE
for the case when the MoC employs the SE solver. In Sections 4 and 5 we repeat
steps of Section 3 for the cases where the ODE solvers are the ME and LF,
respectively. 

It should be noted that in all three cases, the von Neumann analysis leads
one to $4\times 4$ generalized eigenvalue problems which involve both matrices
${\bf B}$ and ${\bf A}$ (or, as per \eqref{e04a}, ${\bf \Sigma}$). 
The largest-magnitude eigenvalue, which determines the stability of the 
numerical method, of those problems cannot be analytically found or 
even related to the {\em individual} 
eigenvalues of the ``physical" matrices ${\bf B}$ and ${\bf A}$. 
However, for given ${\bf B}$ and ${\bf A}$, that eigenvalue can be
easily and quickly found numerically by standard built-in commands in
 software like Matlab, and thus we use this semi-analytical approach.

In Section 6 we comment on the case when the ODE solver is the
fourth-order classical Runge--Kutta (cRK) method. We will not, however, analyze it
in any detail because, as we will demonstrate, there is an uncertainty about
using the cRK solver within the MoC framework. A detailed investigation of this
issue is outside the scope of this paper. In Section 7 we present our
conclusions and demonstrate the validity of our analysis for a broader class
of systems than the particular system considered in Sections 2--5. 
First, in 
Appendix A we apply our analysis of the MoC-LF scheme to a 
class of models of which the model considered in the main text is a special case.
These models still result in linearized equations with constant coefficients.
Then, in Appendix B, we present a different
model with a {\em spatially localized} solution
and demonstrate by direct numerical simulations that essential conclusions 
of our von Neumann analysis remain valid for (in)stability of the MoC applied
to that non-constant-coefficient system.


\section{Physical model}

While our study will focus on a linear problem of a rather general form
(see Eq.~\eqref{e10a} below), we will begin by stating a specific nonlinear
problem which had originally motivated our study and whose linearization
leads to \eqref{e10a}.
We consider the system
\be
\vSpm_{\;t} \pm \vSpm_{\;x} = \vSpm \times {\bf \hat{J}} \vSmp,
\label{e05}
\ee
where $\vSpm\equiv [\Spm_1,\Spm_2,\Spm_3]^T$, \ 
${\bf \hat{J}}={\rm diag}(1,-1,-2)$, and superscript `T' denotes the transposition.
This system is a representative of a class of models that arise in studying
propagation of light in birefringent optical fibers
with Kerr nonlinearity \cite{ZM87}--\cite{Kozlov11}; see Appendix A for more detail. 
The nonlinear system \eqref{e05} has a rather special form and arises in a specific
application. However, its {\em linearized} form (see Eq.~\eqref{e10a} below), for 
which we will analyze the stability of the MoC, is quite general. A wide and diverse
range of physical problems lead to the same equation, as we will explain after
Eqs.~\eqref{e11}.

In the component form, system \eqref{e05} is:
\bsube
\be
(\partial_t+\partial_x)\Sp_1 = \Sp_3 \Sm_2 - 2\Sp_2 \Sm_3, 
\label{e06a}
\ee
\be
(\partial_t+\partial_x)\Sp_2 = 2\Sp_1 \Sm_3 + \Sp_3 \Sm_1, 
\label{e06b}
\ee
\be
(\partial_t+\partial_x)\Sp_3 = -(\Sp_1 \Sm_2 + \Sp_2 \Sm_1), 
\label{e06c}
\ee
\be
(\partial_t-\partial_x)\Sm_1 = \Sm_3 \Sp_2 - 2\Sm_2 \Sp_3, 
\label{e06d}
\ee
\be
(\partial_t-\partial_x)\Sm_2 = 2\Sm_1 \Sp_3 + \Sm_3 \Sp_1, 
\label{e06e}
\ee
\be
(\partial_t-\partial_x)\Sm_3 = (\partial_t+\partial_z)\Sp_3 .
\label{e06f}
\ee
\label{e06}
\esube
It has four families of soliton/kink solutions, a special case of 
one of which is (see, e.g., \cite{Wabnitz09}):
%
%
\be
\Spm_1 = \pm \frac{1}{\sqrt{3}} \,\sech(\sqrt{2} x), \qquad
\Spm_2 = \mp \tanh(\sqrt{2} x),  \qquad
\Spm_3 = \sqrt{2} \,\Spm_1.
\label{e07}
\ee
(The other three families differ from \eqref{e07} by combinations of
signs of the solution's components.) However, we will focus on the analysis
of the (in)stability of the MoC applied to a simpler, {\em constant}
solution:
\be
\Spm_{1,3}=0, \qquad \Spm_2=\pm 1. 
\label{e08}
\ee
This solution corresponds to the asymptotic value of \eqref{e07} as
$x\To -\infty$. 
It can be shown (see, e.g., \cite{Zthesis}) 
that this solution is stable in the sense to be defined below.
However, we have found that when simulated by certain ``flavors" of the MoC, it
can be numerically unstable. If that asymptotic, constant solution is
numerically unstable, then so will be any ``physically interesting" 
non-constant solutions, like \eqref{e07}, possessing the asymptotics
\eqref{e08}. Thus, numerical stability of the constant solution \eqref{e08}
is necessary for the successful performance of the MoC on the soliton/kink
solution \eqref{e07} and related ones. Moreover, 
considering the numerical stability of the MoC simulating the simpler 
solution \eqref{e08}  will allow us to focus on the behavior of the 
numerical method without being distracted by the complexity of the physical 
model. The methodology of our analysis, as well as at least some of our
general conclusions, will hold for the MoC applied to other 
non-dissipative hyperbolic PDEs; we will discuss this in Section 7
and provide details in the Appendices.

For future reference, we linearize Eqs.~\eqref{e06} on the background of
solution \eqref{e08}:
\be
\Spm_j = \Spm_{j\,0}+\tSpm_j, \qquad j=1,2,3;
\label{e09}
\ee
here $\Spm_{j\,0}$ are the components of the exact solution \eqref{e08}
and $\tSpm_j$ are small perturbations. The linearized system \eqref{e06}
reduces to:
\bsube
\be
\tvS_{t} + {\bf \Sigma}\, \tvS_{x} = {\bf P}\, \tvS, 
\label{e10a}
\ee
\be
(\tSpm_2)_t \pm (\tSpm_2)_x = 0,
\label{e10b}
\ee
\label{e10}
\esube
where $\tvS=[\tSp_1,\tSp_3,\tSm_1,\tSm_3]^T$, the $4\times 4$ matrices
in \eqref{e10a} are:
\bsube
\be
{\bf \Sigma} = {\rm diag}(I, -I), \qquad 
{\bf P}= \left( \ba{rr} -A & B \\ -B & A \ea \right), 
\label{e11a}
\ee
and the $2\times 2$ matrices in \eqref{e11a} are:
\be
I = \left( \ba{rr} 1 & 0 \\ 0 & 1 \ea \right), \qquad
A = \left( \ba{rr} 0 & 1 \\ -1 & 0 \ea \right), \qquad
B = -\left( \ba{rr} 0 & 2 \\ 1 & 0 \ea \right).
\label{e11b}
\ee
\label{e11}
\esube

Non-dissipative linearized equations of the form \eqref{e10a}, with 
the same or similar ${\bf \Sigma}$ but a 
more general form of ${\bf P}$, arise in a wide range of physical 
applications involving the interaction of waves propagating with distinctly
different velocities. Examples include: acousto-optical interactions, known
as Stimulated Brillouin Scattering (\cite{Boyd}, Sec.~8.3; 
\cite{SBS_Kaup, SBS_JQE}); \  Stimulated Raman Scattering (\cite{Boyd}, Sec.~9.4;
\cite{SRS_1975}) and, more generally, parametric interaction among two or more
electromagnetic waves \cite{RKB_1979} in semi-classical optics \cite{SGH_Kaup, Ibragimov},
quantum optics \cite{CJM_quant1, CJM_quant2}, photonics \cite{CJM_photon},
and plasma physics applications \cite{CJM_plasma1, CJM_plasma2}; \ 
interaction of crossed beams in photorefractive materials (\cite{Boyd}, Sec.~10.6); \ 
interaction of counter- and co-propagating light waves in a medium with a periodic
refractive index (\cite{FBGbook, 2002_CW, 2014_ChowdAtai}); \ 
and the relativistic 
field theory \cite{Thirring, GrossNeveu} (also see \cite{2006_ChugunovaPeli} 
for a general form of coupled-wave-like field models with a cubic nonlinearity and
\cite{2014_NLDE} for more recent references on these models' solutions).
Interestingly, models similar to those in the relativistic field theory 
also arise in nonlinear fiber optics \cite{AcevesWabnitz, 1992_RTW}.
In Section 7 we will demonstrate that our results for system \eqref{e10},
\eqref{e11} are also applicable to the one-dimensional Gross--Neveu model 
\cite{GrossNeveu} of the relativistic field theory.

In Eqs.~\eqref{e10} and \eqref{e11} and 
in what follows we will adopt the following notations. 
Boldfaced quantities with an underline, ${\bf \vSpm}$, and with a hat,
${\bf \hat{J}}$, will continue to denote 
$3\times 1$ vectors and $3\times 3$ matrices, respectively, as in \eqref{e05}. 
Boldfaced quantities {\em without}
an underline or a hat will denote $4\times 4$ matrices or $4\times 1$ vectors, as
in \eqref{e10a}; the ambiguity of the same notations for matrices and vectors 
here will not
cause any confusion. Finally, underlined letters in regular 
(not boldfaced) font will denote $2\times 1$ vectors; e.g.:
$$
\tvSpm \equiv [\tSpm_1, \tSpm_3]^T.
$$
Clearly then, $\tvS\equiv \left[ (\tvSp)^T,\, (\tvSm)^T \right]^T$.

Seeking the solution of \eqref{e10a} to be proportional to $e^{ikx-i\omega t}$, one
can show that $\omega\in\mathbb{R}$ for all $k\in \mathbb{R}$. This means that
solution \eqref{e08} is stable on the infinite line, as mentioned two paragraphs
above. In particular, $\omega(k=0)\in\mathbb{R}$, which is equivalent to the
statement that all nonzero eigenvalues of ${\bf P}$ are purely imaginary. 
Indeed, from \eqref{e11} one finds:
\be 
\lambda_{\,\bf P} = 0,0,\pm i \sqrt{6}\,; \qquad {\rm so}\quad 
\lambda_{\,\bf P}\in i\,\mathbb{R}.
\label{e12}
\ee

Returning to the full set of equations \eqref{e05} and scalar-multiplying each
of them by its respective $\vSp$ or $\vSm$, we notice that they admit the
following ``conservation" relations:
\bsube
\be
(\partial_t \pm \partial_x )\, |\vSpm|^2 = 0,
\label{e13a}
\ee
where $|\ldots|$ stands for the length of the vector. Two other conservation
relations can be obtained in a similar fashion. Note that for periodic BC,
considered in this paper, these relations become conservation laws. 
E.g., \eqref{e13a} yield:
\be
\partial_t \int_0^L |\vSpm|^2 dx = 0, 
\label{e13b}
\ee
\label{e13}
\esube
where $L$ is the length of the spatial domain. Moreover, a Hamiltonian of system
\eqref{e05} can be constructed in certain action-angle variables \cite{Assemat12}.

Therefore, it is desirable that the numerical method also conserve or quasi-conserve
(i.e., not lead to a systematic shift) at least some of these quantities. A simple,
explicit such a method for ODEs is LF, whereas explicit Euler methods
are known to lead to mild numerical instability for conservative ODEs. For that
reason we had expected that using the LF as the ODE solver in the MoC would produce
better results than explicit Euler solvers. However, we found that, on the contrary, 
it produces the worst results. We will now turn to the description and analysis of
the SE, ME, and LF solvers for the MoC with periodic BC. We will refer to the
corresponding ``flavors" of the MoC as MoC-SE, MoC-ME, and MoC-LF, respectively.


\section{MoC with the SE solver}

The SE is a first-order method and, moreover, is well-known to lead to a mild
yet conspicuous numerical instability when applied to conservative ODEs
(see, e.g., \cite{Griffiths_book}). 
The reason that we consider this method is that it is simple enough to illustrate 
the approach. Thus, describing sufficient details in this Section will 
allow us to skip them 
in subsequent sections, devoted to second-order methods, where such details are
more involved. 


\subsection{Analysis}

The form of the MoC-SE equations for system \eqref{e06} is:
\be
(\Spm_j)^{n+1}_m = (\Spm_j)^{n}_{m\mp 1} + 
      h \, f^{\pm}_j\big(\,  (\vSp)^n_{m\mp 1},\, (\vSm)^n_{m\mp 1}\,\big), 
			\qquad j=1,2,3;
\label{e14}
\ee
where $f_j^{\pm}$ are the nonlinear functions on the rhs of \eqref{e06}, and
$m=0, \ldots\,\, M$, with $(M+1)$ being the number of grid points. Note that in
this paper we have set the temporal and spatial steps equal, $\D t=\D x=h$,
to ensure having a regular grid. 
To impose periodic BC, we make the following identifications in \eqref{e14}:
\be
(\vSpm)^n_{-1} \equiv (\vSpm)^n_M\,, \qquad
(\vSpm)^n_{M+1} \equiv (\vSpm)^n_0\,.
\label{e15}
\ee
Linearizing Eqs.~\eqref{e14}, one arrives at:
\bea
\left( \ba{c} \tvSp \\ \tvSm \ea \right)_m^{n+1} & = &  
\left( \ba{c} \tvSp \\ \und{0} \ea \right)_{m-1}^{n} + 
\left( \ba{c} \und{0} \\ \tvSm \ea \right)_{m+1}^{n} + 
\nonumber \\
& & 
\left( \ba{ll} P^{++} & P^{+-} \\ \mathcal{O} & \mathcal{O} \ea \right) 
\left( \ba{c} \tvSp \\ \tvSm \ea \right)_{m-1}^{n} + 
\left( \ba{ll} \mathcal{O} & \mathcal{O} \\ P^{-+} & P^{--}  \ea \right) 
\left( \ba{c} \tvSp \\ \tvSm \ea \right)_{m+1}^{n},
\label{e16}
\eea 
where $\mathcal{O}$ is the $2\times 2$ zero matrix and
\be
\left( \ba{ll} P^{++} & P^{+-} \\ P^{-+} & P^{--} \ea \right) \equiv 
\left( \ba{ll} -A & B \\ -B & A \ea \right) = 
{\bf P}\,.
\label{e17}
\ee
Seeking the solution of \eqref{e16} with periodic BC in the form
 \ $(\tvS)_m^n = (\tvS)^n\,e^{ikmh}$ \ reduces \eqref{e16} to:
\bsube
\be
(\tvS)^{n+1} = {\bf N}(z)\, (\tvS)^n,
\label{e18a}
\ee
\be
{\bf N}(z) = {\bf Q}\,({\bf I}+h{\bf P}), \qquad 
{\bf Q} = \exp[\,-iz{\bf \Sigma}\,],
\label{e18b}
\ee
\label{e18}
\esube
where $z=kh$, and ${\bf \Sigma}$ is defined in \eqref{e11}.

Below we will show that eigenvalues of ${\bf N}$ satisfy $|\lambda_{\,{\bf N}}|>1$,
which implies that the MoC-SE with periodic BC is unstable. One can easily see this
for the longest- and shortest-period Fourier harmonics
($z=0$ and $z=\pi$, respectively), 
where one has ${\bf Q}(z=0)={\bf I}$ and ${\bf Q}(z=\pi)=-{\bf I}$
and invokes relation \eqref{e12}. Then one has:
\bsube
\be
|\lambda_{\,{\bf N}}(z=0,\pi)| \,=\, \big|1 + ih | \lambda_{\,{\bf P}}|\,\big|\,>1.
\label{e19a}
\ee
However, for other Fourier harmonics of the numerical error, i.e., for
$z\neq 0$ or $\pi$, 
an analytical relation between the eigenvalues
$\lambda_{\,{\bf N}}$ of the numerical method and the 
eigenvalues  $\lambda_{\,{\bf P}}$  of the
``physical" matrix ${\bf P}$ cannot be established. 
In particular, 
\be
|\lambda_{\,{\bf N}}(z\neq 0,\pi)| \;\neq \;|1 + h \lambda_{\,{\bf P}}|,
\label{e19b}
\ee
\label{e19}
\esube
other than by accident. 
Therefore, to determine $|\lambda_{\,{\bf N}}|$, one has to find these
eigenvalues numerically. Let us stress that this limits only the result,
but {\em not the methodology}, of our von Neumann analysis to a specific
matrix ${\bf P}$. Indeed, for any ${\bf P}$, the relation between 
the amplification matrix ${\bf N}$ of the numerical scheme and the
``physical" matrix ${\bf P}$ is given by \eqref{e18b}. The eigenvalues
of ${\bf N}$ are found within a second by modern software,
whereas direct numerical simulations of scheme \eqref{e14}
(and other MoC schemes considered in Sections 4 and 5) take much longer.

\begin{figure}[h!]
\hspace*{-0cm} 
\includegraphics[height=4.4cm,width=5.2cm,angle=0]{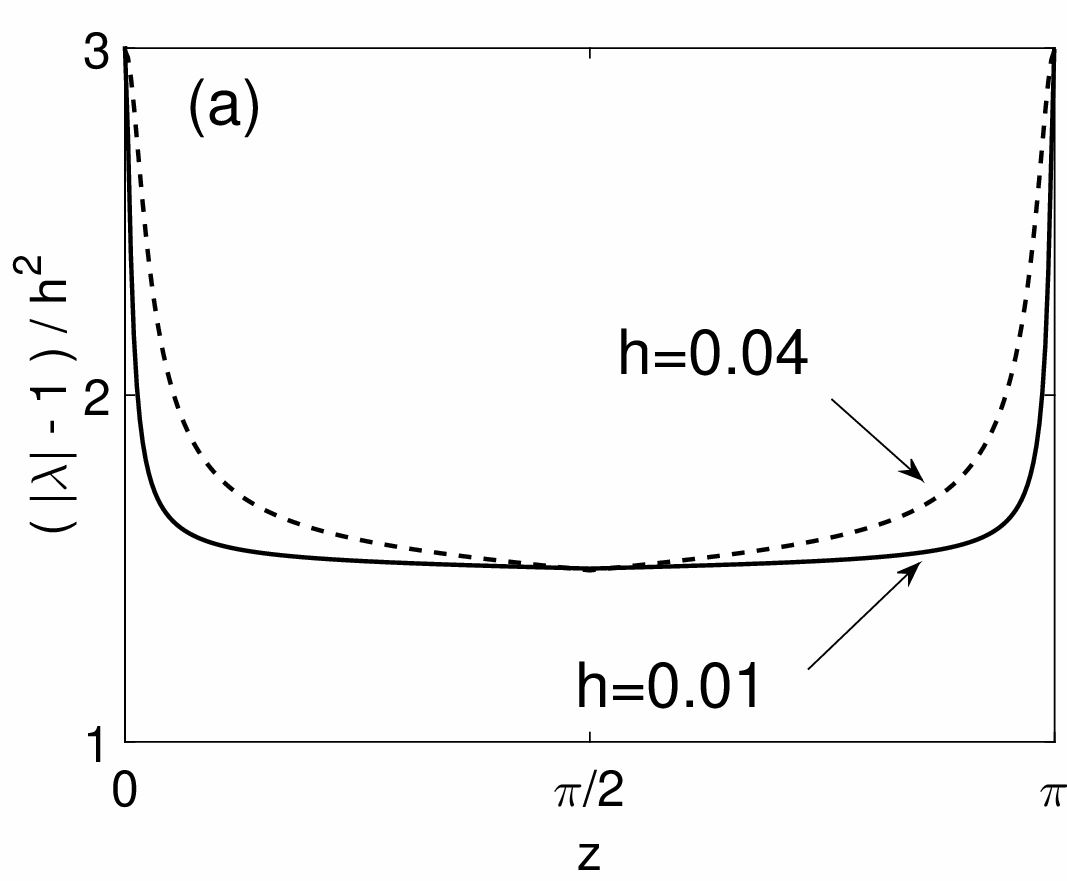}
\hspace{0.1cm}
\includegraphics[height=4.4cm,width=5.2cm,angle=0]{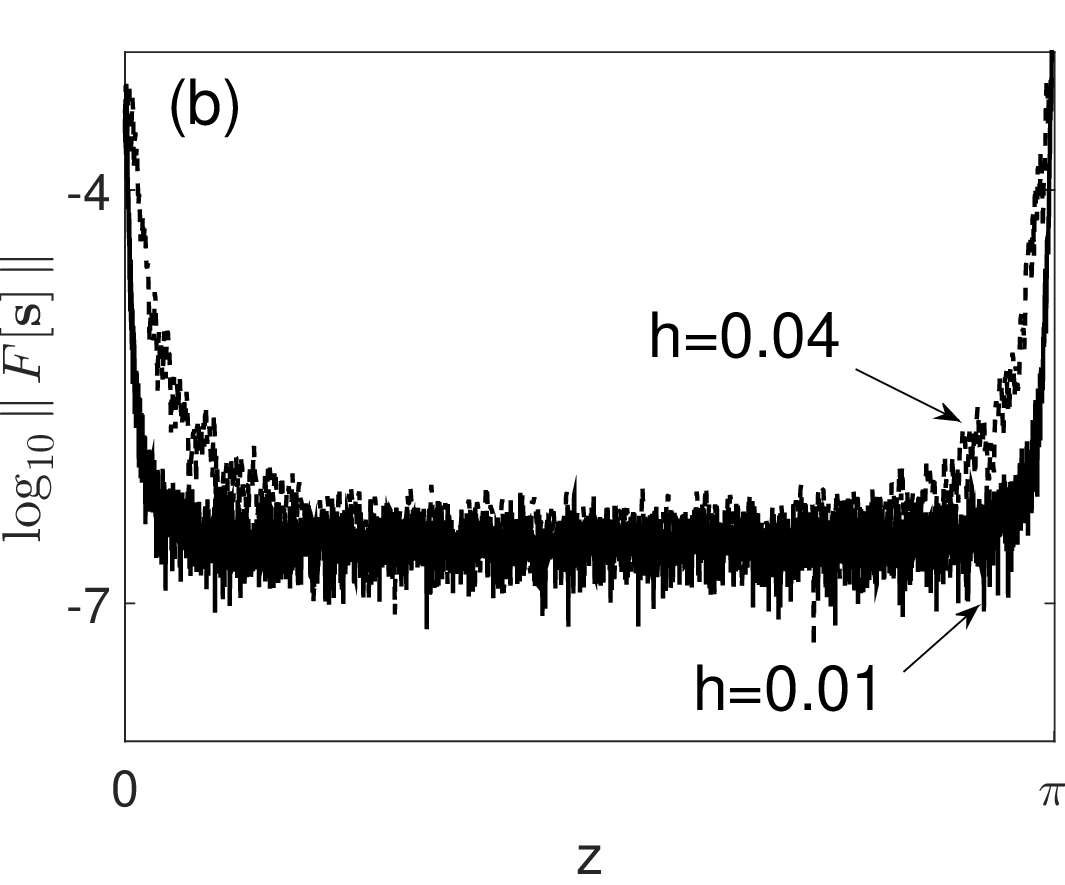}
\hspace{0.1cm}
\includegraphics[height=4.4cm,width=5.2cm,angle=0]{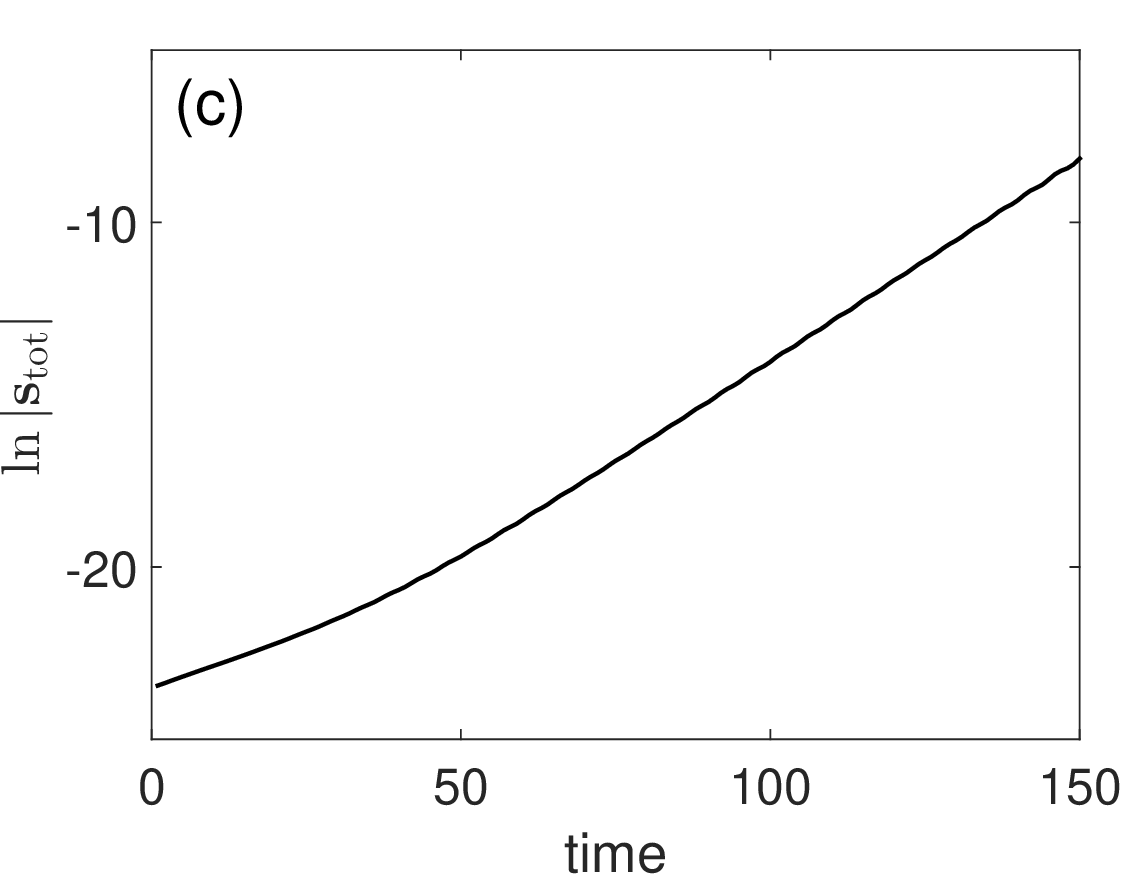}
\caption{
(a) Eigenvalue of matrix ${\bf N}$ in \eqref{e18} for $h=0.01$
and $h=0.04$. \ (b) Spectrum of the numerical error obtained by
simulating scheme \eqref{e14}, \eqref{e15} for $L=100$ and 
$t=520$ (for $h=0.01$) and $t=130$ (for $h=0.04$). \ 
(c) Time evolution of the error \eqref{e20} for $h=0.04$. 
}
\label{fig_1}
\end{figure}

The largest (in magnitude)
eigenvalue of ${\bf N}$, 
found by Matlab, is shown as a function of the normalized wavenumber $z$
in Fig.~\ref{fig_1}(a).
The two curves illustrate the general trend as $h$ is varied. They also show
that the strongest instability of the MoC-SE applied to \eqref{e05} occurs in
the ODE- and ``anti-ODE" limits, i.e. for $z=0$ and $z=\pi$. In particular,
since ${\bf Q}(z=\pi) = -{\bf Q}(z=0)$, one concludes that 
$|\lambda_{\,{\bf N}}(z=\pi)|\,=\,|\lambda_{\,{\bf N}}(z=0)|$. 
In other words, the instability of the MoC-SE for the highest and lowest 
Fourier harmonics is the same, in contrast to that of other numerical methods. 

\subsection{Numerical verification}

First, in Fig.~\ref{fig_1}(b) we show the logarithm of the Fourier spectrum of
the  numerical error obtained by applying the MoC-SE
to \eqref{e06}. Namely, we simulated scheme \eqref{e14} with the initial
condition being \eqref{e08} plus white (in space) noise of magnitude
$10^{-12}$ and subtracted from the numerical solution the exact solution
\eqref{e08}. The curves in Fig.~\ref{fig_1}(b) have the same qualitative
shapes as the corresponding curves in Fig.~\ref{fig_1}(a), which confirms
the validity of the analysis in Section 3.1. Indeed, it follows from
\eqref{e18a} that (individual) Fourier harmonics of the numerical error
satisfy
\be
\| (\tvS)^{n} \| \propto |\lambda_{\,{\bf N}}|^n 
\approx \exp\left[\, \left( (|\lambda_{\,{\bf N}}|-1)/h\right)\,t\,\right], 
\qquad \so \qquad
\ln \| (\tvS)^{n} \| \propto \left( (|\lambda_{\,{\bf N}}|-1)/h\right)\,t\,,
\label{add1_01}
\ee
where $\|\ldots\|$ denotes the Euclidean norm of the vector.

Second, to  verify our analytical results from yet another perspective, 
we will compare
the growth rate measured for a ``total" numerical error (see below)
 with the growth rate predicted by
the analysis of Section 3.1. 
The ``total" numerical error is computed as 
\be
|\tvS_{\rm tot}| = \left( \sum_{m=0}^{M} \| (\tvS)_m \|^2 \right)^{1/2}\,.
\label{e20}
\ee
Note that the summation over $m$ smooths out the noisy spatial profile
of the error. The perturbations $\tSpm_2$ were indeed found to be much
smaller than $\tSpm_{1,3}$ (as long as the latter are themselves 
sufficiently small), as predicted by \eqref{e10b}, and hence they are
not included in the ``total" error \eqref{e20}. 
Now, according to Section 3.1, the highest growth rate occurs
for harmonics with $z=0,\;\pi$. It follows from \eqref{e19a} that
\be
|\lambda_{\,{\bf N}}(z=0)|^n = 
\left( |1+ih\max|\lambda_{\,{\bf P}}| \right)^{t/h} 
\approx \exp\left[ \frac{h}2 \max|\lambda_{\,{\bf P}}|^2 t \right].
\label{e21}
\ee
Thus, the theoretical growth rate is
\be
\gamma_{\rm theor} = 3h,
\label{e22}
\ee
where we have used \eqref{e12}.

On the other hand, the growth rate can be measured from the numerical 
solution:
\be
\gamma_{\rm meas} = 
      \frac{ \ln |\tvS_{\rm tot}(t_2)| - \ln |\tvS_{\rm tot}(t_1)|}{t_2-t_1}\,,
\label{e23}
\ee
where $t_{1,2}$ are some times when the dependence of $\ln |\tvS_{\rm tot}|$
on $t$ appears to be linear (as, e.g., in Fig.~\ref{fig_1}(c) for $t>40$). 
Using \eqref{e23}, we have found for the parameters listed in the caption
to Fig.~\ref{fig_1} that $\gamma_{\rm meas}$ agrees with its theoretical
value \eqref{e22} to two significant figures.


\section{MoC with the ME solver}

\begin{figure}[!ht]
\hspace*{-0cm} 
\includegraphics[height=4.4cm,width=5.2cm,angle=0]{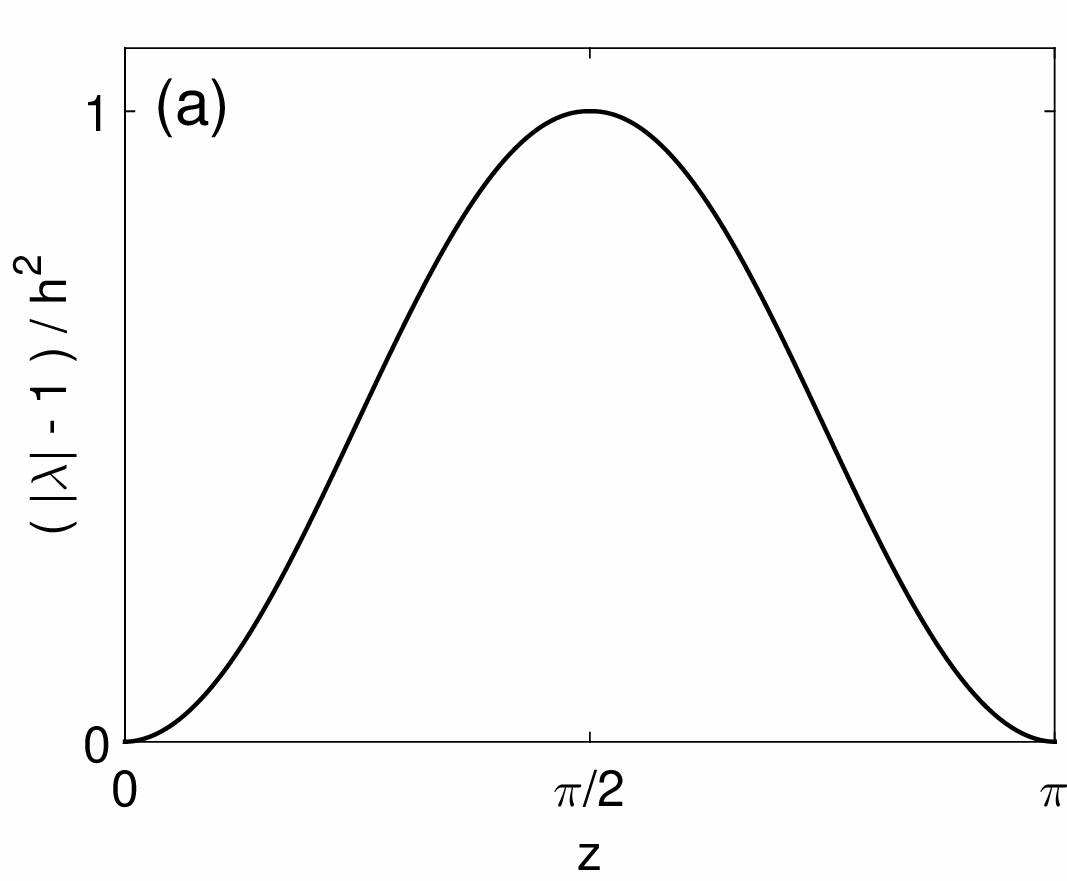}
\hspace{0.1cm}
\includegraphics[height=4.4cm,width=5.2cm,angle=0]{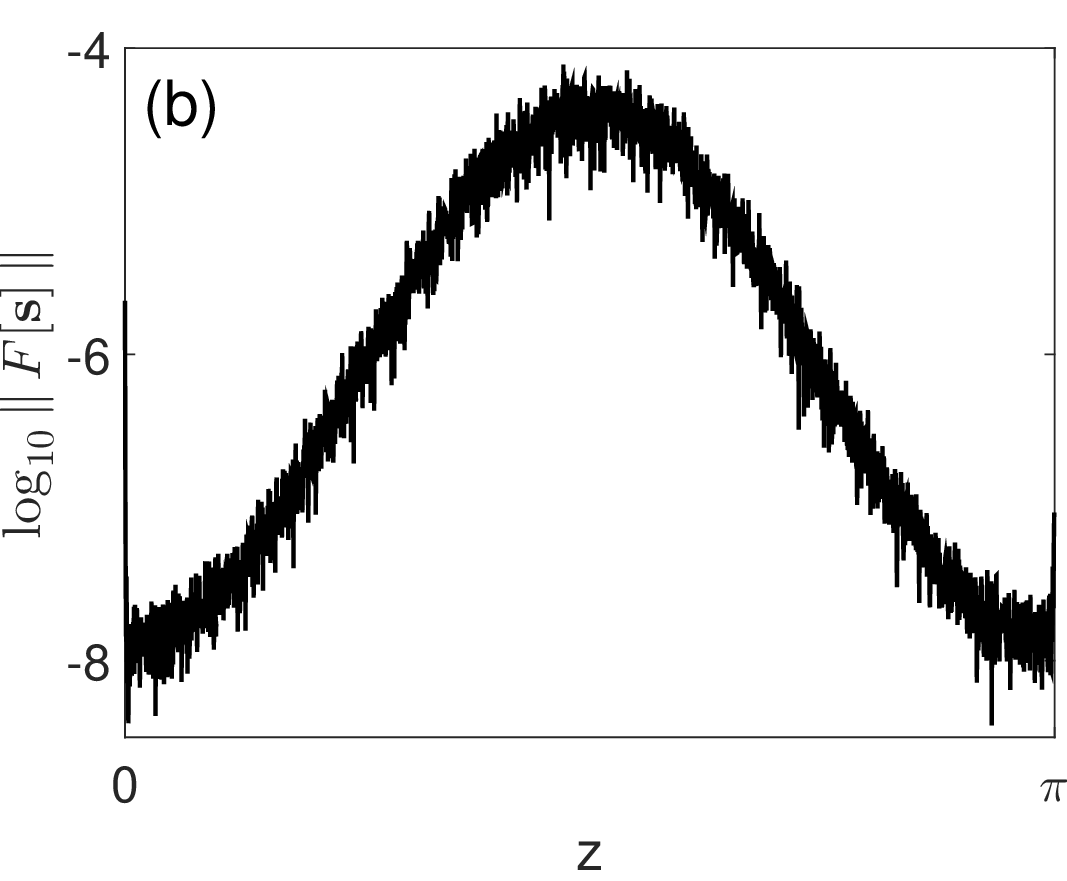}
\hspace{0.1cm}
\includegraphics[height=4.4cm,width=5.2cm,angle=0]{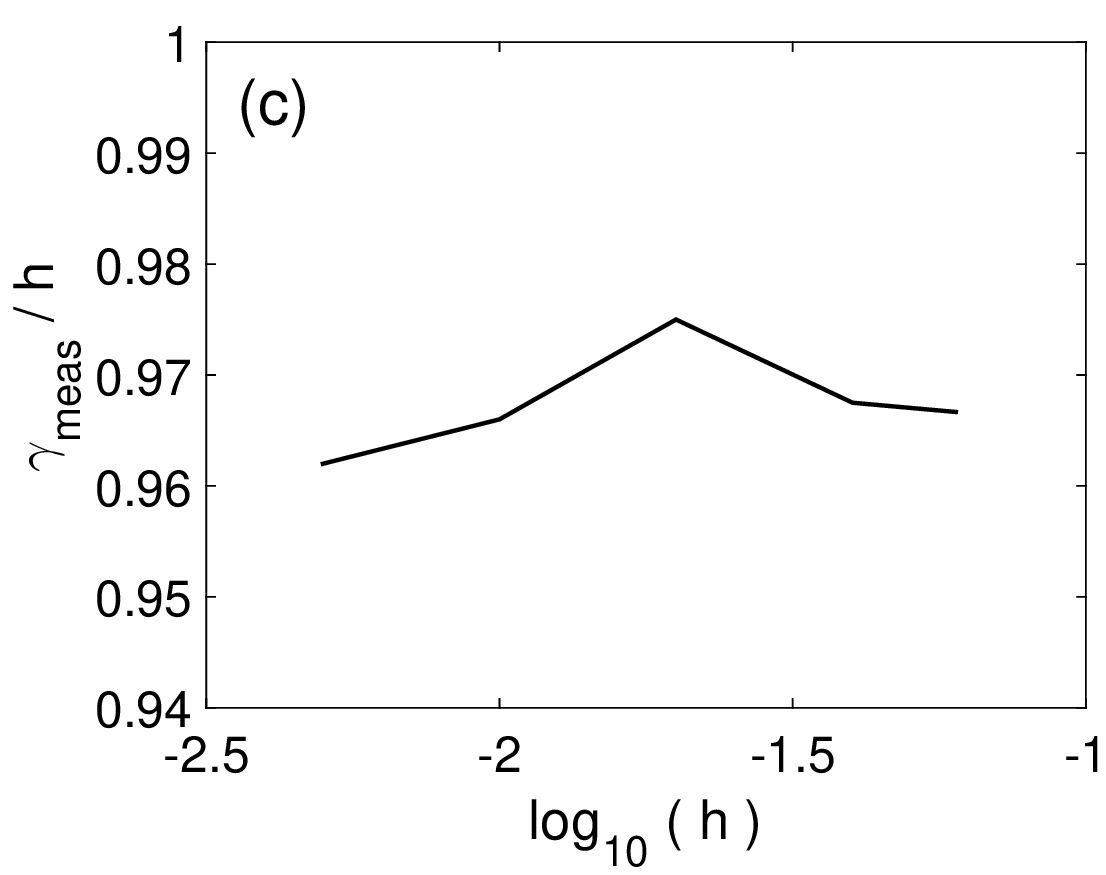}
\caption{
(a) Eigenvalue of matrix ${\bf N}$ in \eqref{e25} for $h=0.01$.
The curves for $h \le 0.05$ all look qualitatively similar to the one shown here;
 in particular, their maximum value is 1. 
 \ (b) Spectrum of the numerical error obtained by
simulating scheme \eqref{e24}, \eqref{e15} for $L=100$ and 
$t=800$ with $h=0.01$. \ 
(c) Comparison of the growth rate of MoC-ME measured from scheme 
\eqref{e24}, \eqref{e15}
with that established by formula \eqref{e26b}. Note that the value 1 of the
ratio would correspond to the perfect agreement between the numerics and
analysis.
}
\label{fig_2}
\end{figure}

The MoC-ME applied to system \eqref{e06} yields the following scheme:
\bsube
\be
\overline{\Spm_j}_m  =   (\Spm_j)^{n}_{m\mp 1} + 
      h \, f^{\pm}_j\big(\,  (\vSp)^n_{m\mp 1},\, (\vSm)^n_{m\mp 1}\,\big), 
			\qquad j=1,2,3;
\label{e24a}
\ee
\be
(\Spm_j)^{n+1}_m  =   \frac12 \left[ \, 
                      (\Spm_j)^{n}_{m\mp 1} + \overline{\Spm_j}_m + 
      h \, f^{\pm}_j\Big(\,  \overline{\vSp}_m,\, \overline{\vSm}_m\,\Big)
			 \, \right]\,, 
\label{e24b}
\ee
\label{e24}
\esube
where the notations are the same as in Section 3. Note that the rhs of 
\eqref{e24a} is the same as that of \eqref{e14}. Following the analysis of 
Section 3.1, one obtains a counterpart of \eqref{e18} where now
\be
{\bf N} = \frac12 \left[ {\bf Q} + 
          ({\bf I}+h{\bf P})\,{\bf Q}\, ({\bf I} + h{\bf P})\, \right].
\label{e25}
\ee
As we discussed in Section 3.1, the eigenvalues of ${\bf N}$ 
cannot be analytically related to the eigenvalues of ${\bf P}$ and hence
need to be found numerically. (As for the MoC-SE before, 
this takes fractions of a second for modern software, unlike the much longer 
direct simulations of scheme \eqref{e24}.)
The magnitude of the largest-in-magnitude such eigenvalue 
is shown in Fig.~\ref{fig_2}(a).
The logarithm of the Fourier spectrum of $\|\tvS\|$ is plotted in 
Fig.~\ref{fig_2}(b) and is seen to agree qualitatively with the result
in Fig.~\ref{fig_2}(a); see also \eqref{add1_01}. 

Let us note that the numerical instability of the MoC-ME scheme in the
ODE and anti-ODE limits is much weaker than that of the MoC-SE, because
from \eqref{e25}, 
$$
|\lambda_{\,{\bf N}}(z=0,\,\pi)| = \left( 
                           1 + \frac{h^4}{4}\max|\lambda_{\,{\bf P}}|^4 \right)^{1/2},
$$
and hence similarly to \eqref{e21}, one has
\bsube
\be
\gamma_{\rm theor}(z=0,\,\pi) = \frac{h^3}{8}\max|\lambda_{\,{\bf P}}|^4\,.
\label{e26a}
\ee
However, Fig.~\ref{fig_2}(a) shows that {\em unlike} the MoC-SE, the
strongest instability of the MoC-ME occurs in the ``{\em middle}"\footnote{
	We used the quotes because, strictly speaking, this is the middle of the
	right {\em half} of the spectrum.
	}
 of the Fourier
spectrum. Moreover, this growth rate has the same 
order of magnitude 
as the growth rate of the MoC-SE's instability.
Therefore, it is the growth rate of the harmonics with 
$k\sim \pi/(2h)$ that one would measure in the numerical experiment. 
Using the information from Fig.~\ref{fig_2}(a) and its caption, one has
$$
|\lambda_{\,{\bf N}}(z=\pi/2)| \approx   1 + h^2,
$$
and hence
\be
\gamma_{\rm theor}(z=\pi/2) \approx h.
\label{e26b}
\ee
\label{e26}
\esube
 Figure \ref{fig_2}(c) shows that the 
growth rates computed from the (semi-)analytical formula
\eqref{e26b} and measured, as explained in Section 3.2, 
in direct numerical simulations of scheme \eqref{e24}, agree
reasonably well.

In Section 7 we will demonstrate that these conclusions
about the numerical instability of the MoC-ME remain qualitatively
valid for a different system of PDEs, which has spatially localized 
(as opposed to constant) coefficients; see Appendix B.


\section{MoC with the LF solver}

\begin{figure}[h!]
\begin{minipage}{7.5cm}
\hspace*{-0.1cm} 
\includegraphics[height=6cm,width=7.5cm,angle=0]{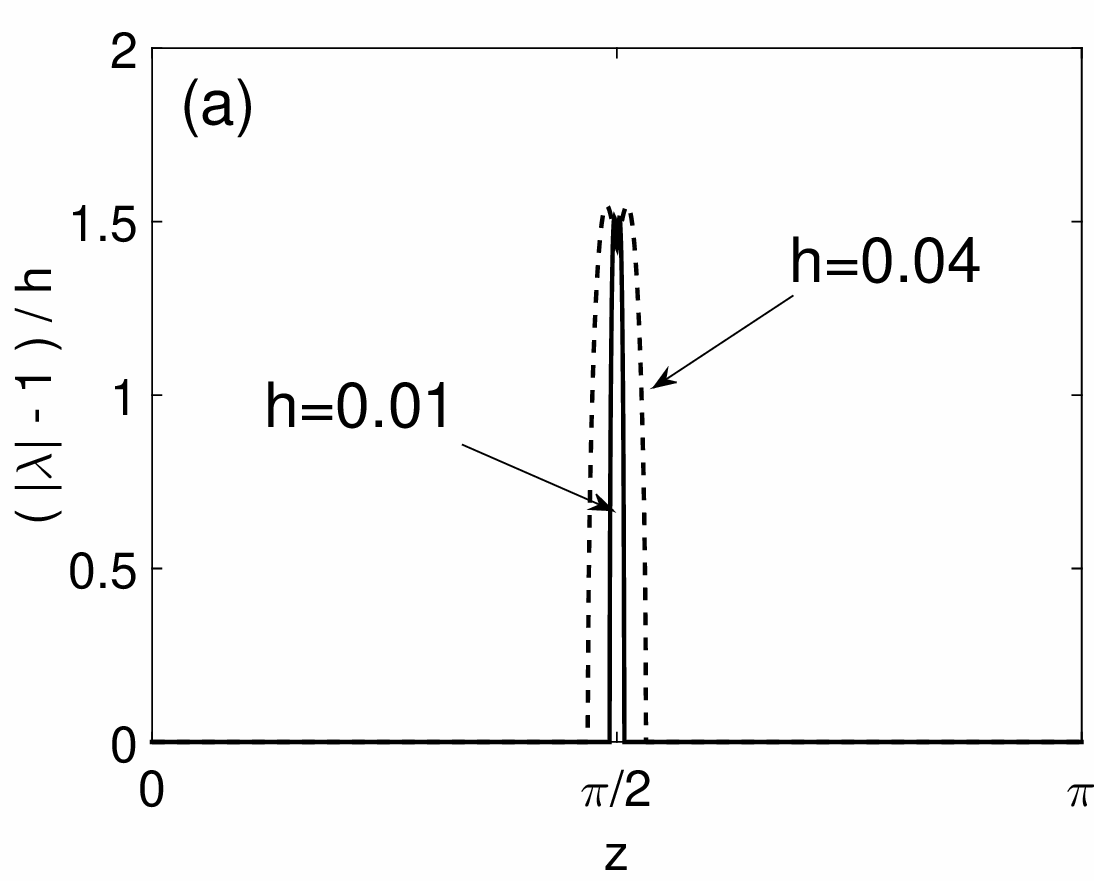}
\end{minipage}
\hspace{0.5cm}
\begin{minipage}{7.5cm}
\hspace*{-0.1cm} 
\includegraphics[height=6cm,width=7.5cm,angle=0]{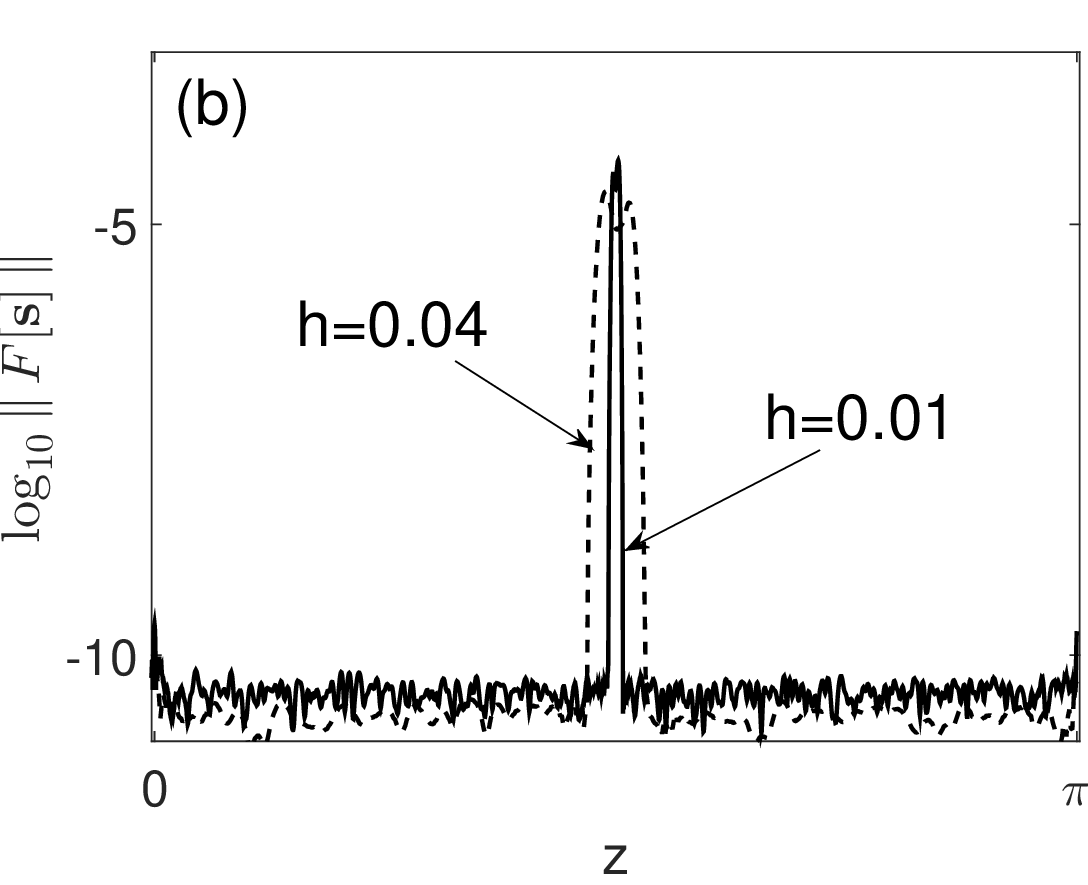}
\end{minipage}
\caption{(a) Eigenvalue of problem \eqref{e30} for $h=0.01$
and $h=0.04$. \ (b) Spectrum of the numerical error obtained by
simulating scheme \eqref{e27}, \eqref{e15}, \eqref{e28} for $L=100$ and 
$t=10$. \ 
}
\label{fig_3}
\end{figure}

Unlike the Euler ODE solvers considered in the previous two sections,
the LF involves the solution at three time levels:
\be
(\Spm_j)^{n+1}_m = (\Spm_j)^{n}_{m\mp 2} + 
      2h \, f^{\pm}_j\big(\,  (\vSp)^n_{m\mp 1},\, (\vSm)^n_{m\mp 1}\,\big), 
			\qquad j=1,2,3.
\label{e27}
\ee
Given the initial condition, one can find the solution at time level $n=1$
with, e.g., the MoC-SE. To enforce periodic BC, convention \eqref{e15} needs
to be supplemented by
\be
(\vSpm)^n_{-2} \equiv (\vSpm)^n_{M-1}\,, \qquad
(\vSpm)^n_{M+2} \equiv (\vSpm)^n_1\,.
\label{e28}
\ee
The counterpart of \eqref{e18} is now
\be
(\tvS)^{n+1} - 2h{\bf Q}{\bf P}\,(\tvS)^n + {\bf Q}^2 \,(\tvS)^{n-1} = {\bf 0}.
\label{e29}
\ee
The (in)stability of the MoC-LF is determined by whether the magnitude
of the largest eigenvalue of the quadratic eigenvalue problem obtained
from \eqref{e29} and leading to 
\be
\det \left( \lambda^2{\bf I} - 2\lambda h {\bf Q\,P} + {\bf Q}^2 \right) = 0
\label{e30}
\ee
exceeds unity. 
As previously for the MoC-SE and MoC-ME, this $4\times 4$ eigenvalue problem
has to be solved numerically for a given ${\bf P}$.
The largest magnitude of its eigenvalue, computed with
Matlab's command \verb+polyeig+ for ${\bf P}$ given by \eqref{e11}, 
is plotted in Fig.~\ref{fig_3}(a).

It should be noted that the quantity $(|\lambda|-1)$ for the MoC-LF 
scales as $h$, in contrast to the cases of MoC-SE and MoC-ME, where it
scales as $h^2$. Therefore, the instability growth rate of the MoC-LF
is: $\gamma=O(1)$, which is much greater than those rates of the
MoC-SE and MoC-ME. The reason for this will come out as a by-product
when we qualitatively explain, in the next paragraph, why the 
strongest instability of the MoC-LF for the PDE system in question
occurs near $z=\pi/2$.

First note that in the ODE limit, $z=0$, one has ${\bf Q}(0)={\bf I}$,
and hence from \eqref{e30} one recovers the well-known relation between
the amplification factor $\lambda$ of the ODE-LF and eigenvalues of ${\bf P}$:
\be
\lambda(z=0) = h\lambda_{\,{\bf P}} \pm 
               \sqrt{ 1 + (h\lambda_{\,{\bf P}})^2 }.
\label{e31}
\ee
Given $\lambda_{\,{\bf P}}\in i\mathbb{R}$, as would be the case for any
stable non-dissipative system, this yields $|\lambda(z=0)|=1$ 
(here and below we assume that $h\max|\lambda_{\,{\bf P}}| < 1$). 
A similar situation holds in the anti-ODE limit, where ${\bf Q}=-{\bf I}$.
However, for $z=\pi/2$, ${\bf Q}(\pi/2)=-i\,{\rm diag}(I,-I)$ can be said
to be the most dissimilar from ${\bf I}\equiv {\rm diag}(I,I)$, 
and the counterpart 
of \eqref{e31} is:
\be
\lambda(z=\pi/2) = -h\lambda_{\,{\bf P}_{\rm mod}} \pm 
               \sqrt{ -1 + (h\lambda_{\,{\bf P}_{\rm mod}})^2 },
\label{e32}
\ee
where ${\bf P}_{\rm mod}={\bf Q}(\pi/2)\,{\bf P}$. For 
$\lambda_{\,{\bf P}_{\rm mod}}\in i\mathbb{R}$, this yields
$\max|\lambda(z=\pi/2)|>1$, i.e., a numerical instability. In the case
of matrix ${\bf P}$ defined in \eqref{e11}, 
$\lambda_{\,{\bf P}_{\rm mod}}=0,0,\pm i \sqrt{2}$, and thus the above statements
explain our result shown in Fig.~\ref{fig_3}(a). In Section 7 we will discuss
how this generalizes to some other energy-preserving hyperbolic PDEs.

Moreover, it also follows from 
\eqref{e32} that the growth rate of this instability is about
$\max|\lambda_{\,{\bf P}_{\rm mod}}|=O(1)$. 
(For the specific ${\bf P}$ in \eqref{e11}, we will show in the 
second paper of this work that the growth rate is $3/2$ instead of
$\sqrt{2}=\max|\lambda_{\,{\bf P}_{\rm mod}}|$; this minor difference
is responsible for the subtly double-peaked
 profile of the curves in Fig.~\ref{fig_3}(a).)
We confirmed this by direct
numerical simulations: representative spectra of the numerical error are
shown in Fig.~\ref{fig_3}(b)
and are seen to agree with the shapes of the curves in Fig.~\ref{fig_3}(a);
see \eqref{add1_01}.
 The solution obtained by the MoC-LF becomes
completely destroyed by the numerical instability around $t=25$ regardless
of $h$. In contrast, for the MoC-ME, the solution remains essentially unaffected
by the growing numerical error over times $O(1/h)\gg 10$.


\section{Comments about the MoC with the cRK solver}

The fourth-order cRK method may appear as a natural choice of a higher-order
ODE solver to be combined with the MoC. Indeed, not only is its accuracy for
ODEs $O(\D t^4)$, but also the systematic error that it introduces into the
numerically computed conserved quantities of energy-preserving ODEs is usually
negligible for sufficiently small time steps.

However, a straightforward use of the cRK solver in the MoC framework produces
an unsatisfactory numerical method. The corresponding scheme is:
\bsube
\be
(\vSpm)_{m}^{n+1} = (\vSpm)_{m\mp 1}^{n} + \frac16
              \left( (\und{\boldsymbol{\kappa}}_{\,1}^{\pm})_m + 
					 2(\und{\boldsymbol{\kappa}}_{\,2}^{\pm})_m
	              + 2(\und{\boldsymbol{\kappa}}_{\,3}^{\pm})_m + 
					 (\und{\boldsymbol{\kappa}}_{\,4}^{\pm})_m \right);
\label{e33a}
\ee
\be
(\und{\boldsymbol{\kappa}}_{\,1}^{\pm})_m = h\,\und{\bf f}^{\pm}
  \left(\,\vSp_{m\mp 1}, \; \vSm_{m\mp1}\,\right);
\label{e33b}
\ee
\be
(\und{\boldsymbol{\kappa}}_{\;\alpha}^{\pm})_m = h\,\und{\bf f}^{\pm}
  \left(\,\vSp_{m\mp 1} + \frac12 (\und{\boldsymbol{\kappa}}_{\;\alpha-1}^{+})_m, 
	     \;   \vSm_{m\mp1} + \frac12 (\und{\boldsymbol{\kappa}}_{\;\alpha-1}^{-})_m
					\,\right), 
					\qquad \alpha=2,3;
\label{e33c}
\ee
\be
(\und{\boldsymbol{\kappa}}_{\,4}^{\pm})_m = h\,\und{\bf f}^{\pm}
  \left(\,\vSp_{m\mp 1} +  (\und{\boldsymbol{\kappa}}_{\,3}^{+})_m, 
	     \;   \vSm_{m\mp1} + (\und{\boldsymbol{\kappa}}_{\,3}^{-})_m
					\,\right).
\label{e33d}
\ee
\label{e33}
\esube
A von Neumann analysis analogous to that presented in Section 4 yields a
dependence $\max|\lambda_{\,{\bf N}}(z)|$ similar to the one shown in
Fig.~\ref{fig_2}(a), and direct numerical simulations confirm that result
for system \eqref{e06}.
Thus, scheme \eqref{e33} has the instability growth rate $\gamma=O(h)$,
just as the MoC-ME \eqref{e24}. In fact, for the specific PDE system \eqref{e06},
this $\gamma$ is twice that of the MoC-ME's, and hence scheme \eqref{e33}
offers no advantage over the MoC-ME.

Perhaps, a proper MoC-cRK needs to ensure that 
the arguments of $\und{\bf f}^{\pm}$ in \eqref{e33c} and \eqref{e33d}
be taken on the same characteristics. E.g., \eqref{e33d} would then 
need to be replaced with
$$
(\und{\boldsymbol{\kappa}}_{\,4}^{\pm})_m = h\,\und{\bf f}^{\pm}
  \left(\,\vSp_{m - 1} +  (\und{\boldsymbol{\kappa}}_{\,3}^{+})_m, 
	       \; \vSm_{m +1} + (\und{\boldsymbol{\kappa}}_{\,3}^{-})_m
					\,\right),
$$
which mimics the last term in \eqref{e24}. The node numbers of some
of the arguments of $\und{\bf f}^{\pm}$ in \eqref{e33c} would then
also need to be modified in some way. This may result in a more stable
MoC-cRK. However, a detailed exploration of this topic will have to begin
with an analysis of the accuracy of the resulting method away from the 
ODE limit. This is outside the scope of this paper, and therefore we
will not pursue it here. Let us mention that we have been unable to find
any published research which would study the problem of setting up the
proper MoC-cRK for the case of crossing characteristics.

\section{Summary and discussion}

We have presented results of the von Neumann analysis of the MoC combined
with three ODE solvers: SE, ME, and LF, as applied to an energy-preserving 
hyperbolic PDE system, \eqref{e05} or \eqref{e06}. Our main findings are
as follows.

First, we have found that the numerical instability of the highest and lowest
Fourier harmonics in the MoC schemes has the same growth rate. This is in
contrast with the situation of other numerical methods, where it is the
highest Fourier harmonics that become numerically unstable before the lower ones.
Moreover, we have found that the harmonics in the ``middle" of the spectrum,
i.e., with $|k|\approx \pi/(2h)$, can become the most numerically unstable. 
This occurs, for the PDE system in question, for the MoC-ME and MoC-LF schemes.
While such harmonics for the MoC-LF occupy a narrow spectral band and hence
can, in principle, be filtered out, for the MoC-ME they occupy a substantial 
portion of the spectrum and hence cannot be filtered out without destroying
the solution. It should also be pointed out that at least for some energy-preserving
{\em PDEs}, as the one considered here, 
the growth rate of the most unstable harmonics
has the same order of magnitude, $O(h)$, for the MoC-SE and MoC-ME. This should
be contrasted with the situation for energy-preserving {\em ODEs}, 
where the instability growth rate of the ME method, $O(h^3)$, is much lower
than that, $O(h)$, of the SE method.

Second, we have found that the MoC based on the LF solver, which (the solver)
is known to considerably outperform the Euler schemes for energy-preserving
ODEs, performs the worst for the energy-preserving PDE that we considered.
Namely, the Fourier harmonics with $|k|\approx \pi/(2h)$ grow at the rate
$\gamma=O(1)$, which by far exceeds the growth rate $\gamma=O(h)$
of the numerically most unstable harmonics in both the MoC-SE and MoC-ME.

In Section 6 we have pointed out that,
for hyperbolic PDE systems with crossing characteristics, the formulation
of a MoC based on the cRK ODE solver and having
 numerical stability (or, at worst, a milder instability than
the MoC-ME) for $|kh|\sim \pi/2$, remains an open problem.

There remains a question as to whether the above conclusions of our analysis,
which necessarily had to be obtained for a {\em specific} system (as we explained
in the main text), apply to other energy-preserving hyperbolic PDEs.
Again, for the same reasons, this cannot be answered in general, as no
analytical relation between the eigenvalues of the ``physical" matrix ${\bf P}$
and the numerical amplification factor $|\lambda_{\bf N}|$ can be obtained
for an arbitrary wavenumber $k$. Therefore, in Appendices A and B we consider
several specific PDEs with the same linearized form \eqref{e10a}, where matrix
${\bf P}$ has a different form than \eqref{e11}. In Appendix A we limit ourselves
to systems with spatially constant coefficients and focus on how 
qualitatively different eigenvalues (see below) of ${\bf P}$ affect the
amplification factor for the MoC-LF. In Appendix B we will demonstrate that 
our conclusions about numerical instability of both MoC-LF and MoC-ME
remain qualitatively valid for a different system whose linearization
has spatially-{\em localized} entries in matrix ${\bf P}$.


Thus, we will now  explore what ``patterns" of numerical instability
can be expected in simulations of energy-preserving PDEs by the MoC-LF.
To that end,
we consider a number of systems\footnote{
	By a ``system" we will refer here to the PDE plus the stationary solution
	on whose background we perform the linearization. Both of them are needed to
	define entries in matrix ${\bf P}$.
	}
representing a broader class of such PDEs, which
include our system \eqref{e06}, \eqref{e08} as a special case.
The total of 18 systems, including 
it, were considered, and a summary of results is presented in Figs.~\ref{fig_3}(a)
and \ref{fig_4} and explained below. The description of the systems is found
in Appendix A.

\begin{figure}[!ht]
\hspace*{-0cm} 
\includegraphics[height=4.4cm,width=5.2cm,angle=0]{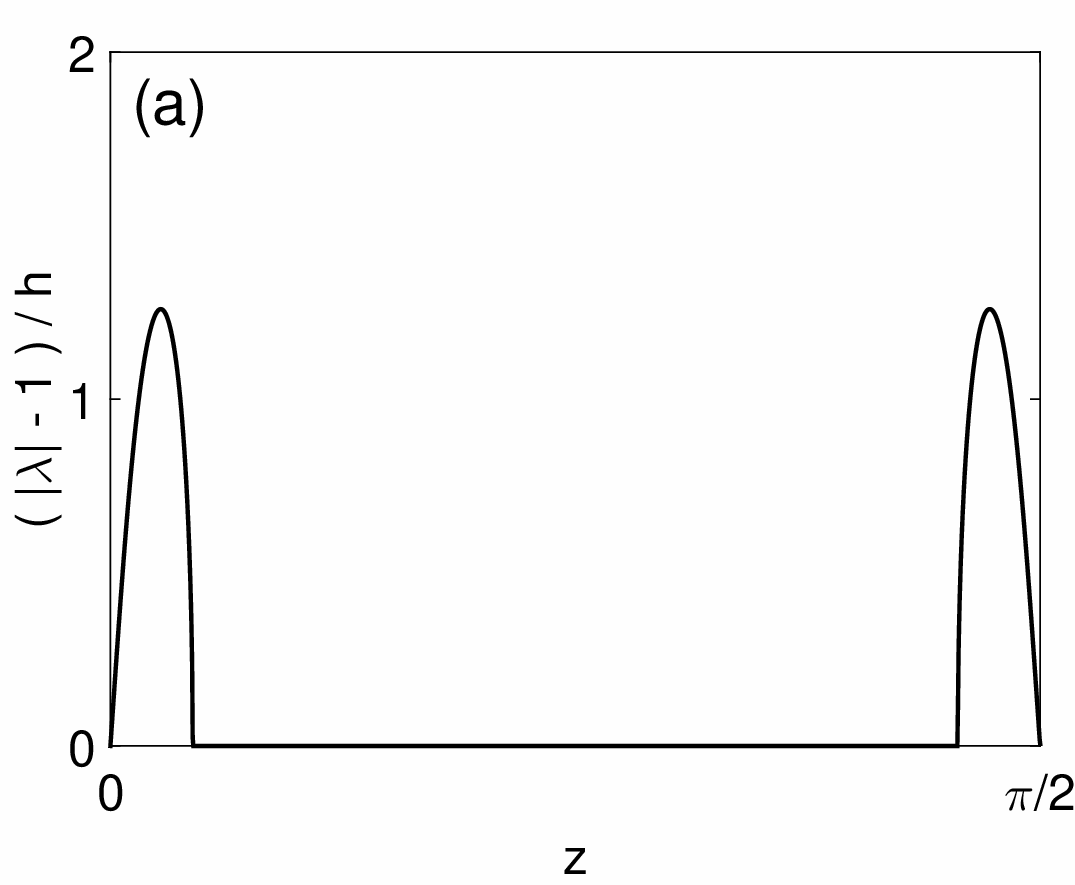}
\hspace{0.1cm}
\includegraphics[height=4.4cm,width=5.2cm,angle=0]{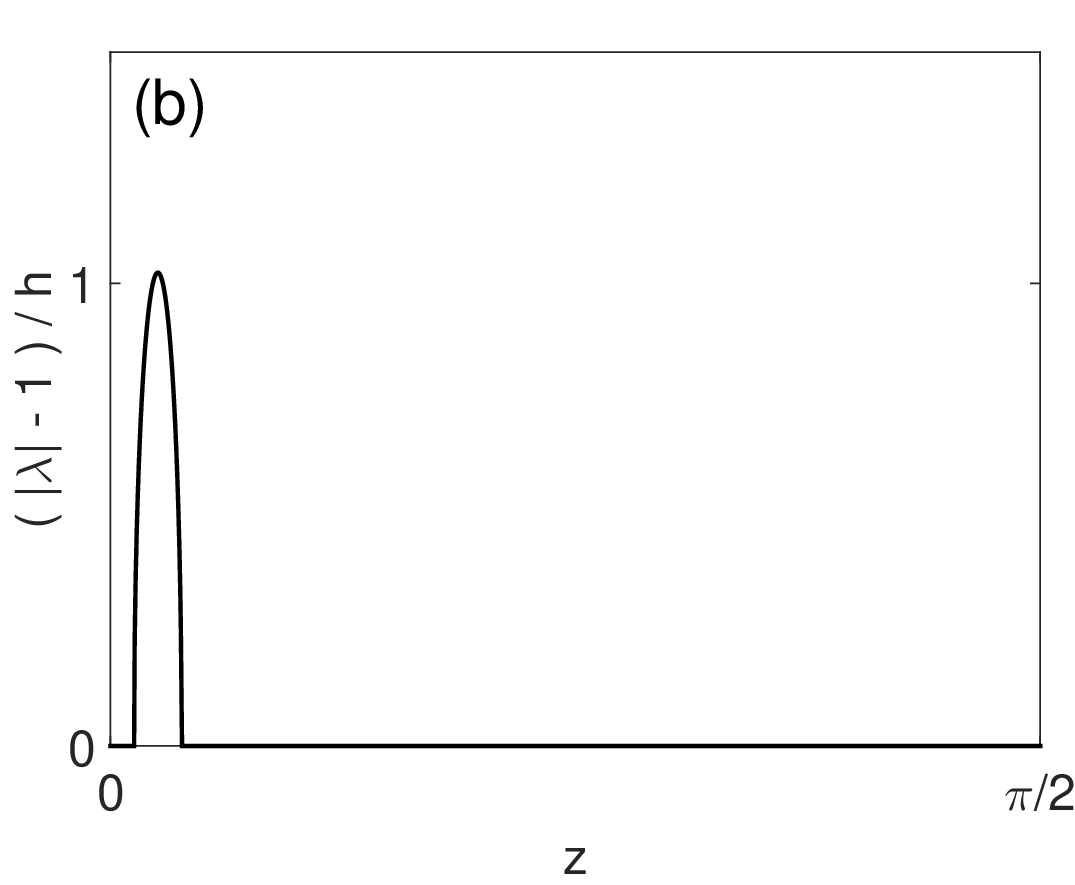}
\hspace{0.1cm}
\includegraphics[height=4.4cm,width=5.2cm,angle=0]{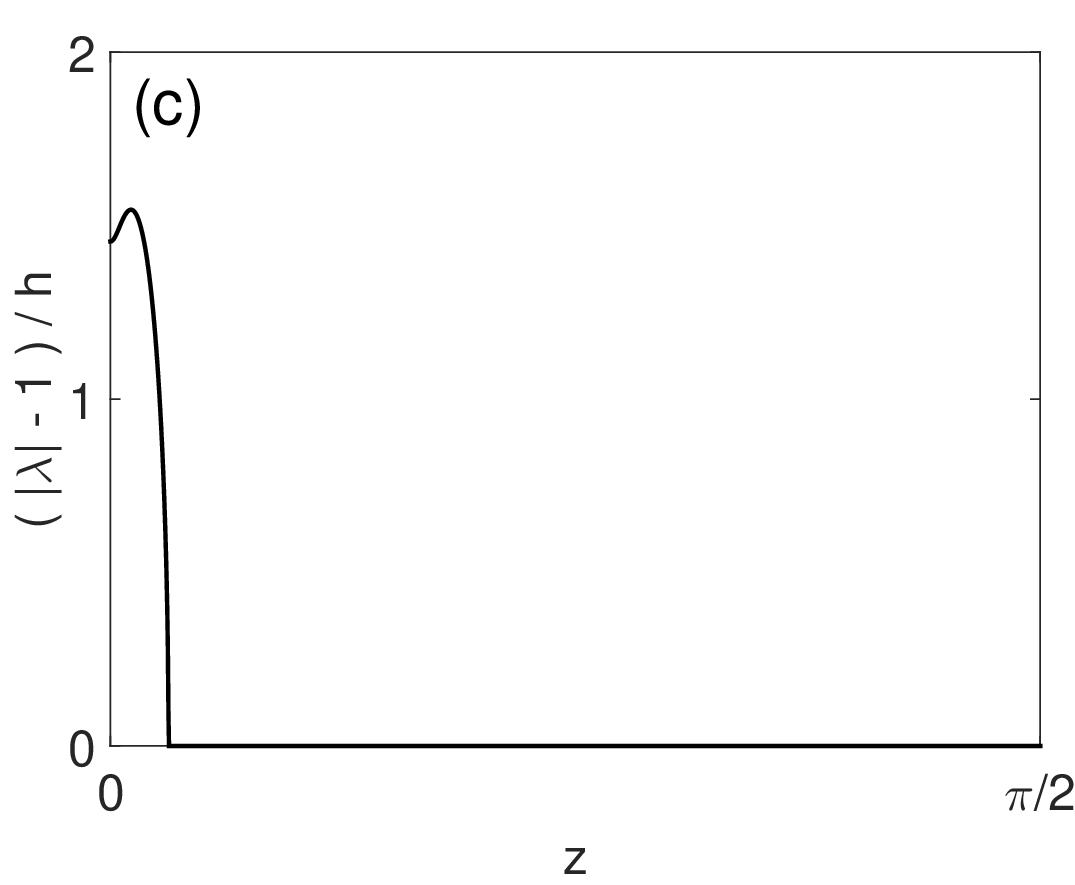}
\caption{
Eigenvalue of problem \eqref{e30} for $h=0.04$ and for the matrices
${\bf P}$ corresponding to the physically 
unstable systems considered in Appendix A. 
The ``other half" of the spectrum, $z\in[\pi/2,\,\pi]$, is reflectionally
symmetric about the point $z=\pi/2$ to the one shown here. 
(a) $\lambda_{\bf P}=\pm i\alpha,\; \pm i\beta$; \ 
 \ (b) $\lambda_{\bf P}=\pm i\alpha$ (a pair of double roots); \ 
(c) $\lambda_{\bf P}=\pm \alpha,\; \pm \beta$. See 
Section 7 for more details. 
}
\label{fig_4}
\end{figure}

Before we proceed, let us remind the reader that the stability of the 
{\em physical} problem \eqref{e10a}, as opposed to the stability of the 
numerical method, was defined before Eq.~\eqref{e12}. For example, system
\eqref{e06}, \eqref{e08} is stable, which can be seen from Fig.~\ref{fig_3}(a).
Indeed, the MoC-LF accurately approximates the solution for $k=O(1)$
(i.e., for $z\equiv kh = O(h)$ in Fig.~\ref{fig_3}(a)), and the fact that
$|\lambda|=1$ there corresponds to $\omega \in \mathbb{R}$ in the 
text before Eq.~\eqref{e12}. Let us stress that our focus is on
the presence of numerical instability of the MoC-LF for physically
{\em stable} systems. Indeed, even if the numerical method is found to 
be free of numerical instability for a physically unstable system, then
the corresponding solution will eventually (typically within $t=O(1)$) 
evolve towards the physically stable solution.

We now continue with the summary of stability results for the MoC-LF.
These were obtained by solving the quadratic eigenvalue problem \eqref{e30}
with Matlab's command \verb+polyeig+ for 18 specific matrices ${\bf P}$
corresponding to each of the systems listed in Appendix A.
For all the seven of those systems which are physically stable,
the dependence $\max |\lambda(z)|$ was found to be similar to that
shown in Fig.~\ref{fig_3}(a). Thus, for all of the stable systems considered,
the MoC-LF has numerically unstable modes with $|kh|\approx \pi/2$ and
growth rate $\gamma=O(1)$. Thus, this method should be deemed unsuitable
for simulation of these, and possibly other, stable energy-preserving 
PDEs with crossing characteristics.

Let us mention, for completeness, that for the 
physically unstable systems considered 
in Appendix A, we found two possibilities with respect to the numerical 
(in)stability of the MoC-LF. In the first group, there are unstable systems
such that: \ $\lambda_{\bf P}=\pm i\alpha,\;\pm i\beta$ for 
$\alpha,\beta\in \mathbb{R}$ (i.e., the mode with $k=0$ is stable),
with $\alpha\neq \beta$ and at most one of $\alpha$ or $\beta$ could be zero.
We have found that such systems have $\omega\,\cancel{\in}\,\mathbb{R}$ in 
the immediate vicinity of $k=0$. For systems in this group, the MoC-LF 
has also numerically unstable modes with $|k|h$ near (but not exactly at) 
$\pi/2$. This is illustrated in Fig.~\ref{fig_4}(a). In the second group
there are unstable systems with all other possibilities of $\lambda_{\bf P}$
compatible with the energy-preserving nature of the problem: \ 
$\lambda_{\bf P}=\pm i\alpha$ (two pairs of repeated imaginary eigenvalues),
including the case $\alpha=0$; and $\lambda_{\bf P}=\pm\alpha,\;\pm \beta$
(in this case, even the $k=0$ mode is unstable). This is illustrated in
Figs.~\ref{fig_4}(b,c), which show {\em no numerical} instability
in the MoC-LF. 
 However, as we have mentioned above, once the initial, physically
unstable, solution evolves sufficiently near a stable one, the MoC-LF will 
be invalidated by the numerically unstable harmonics, which are
seen in Figs.~\ref{fig_3}(a) and  \ref{fig_4}(a).

\begin{figure}[!ht]
\hspace*{-0cm} 
\includegraphics[height=4.4cm,width=5.2cm,angle=0]{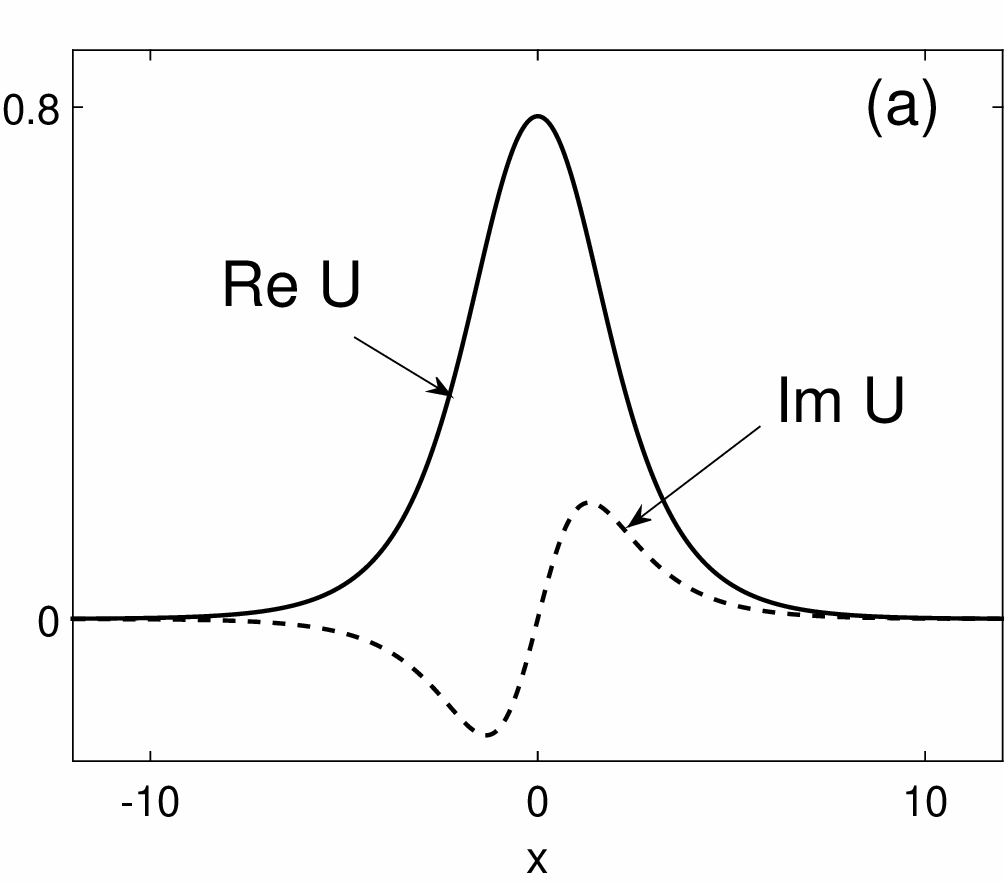}
\hspace{0.1cm}
\includegraphics[height=4.4cm,width=5.2cm,angle=0]{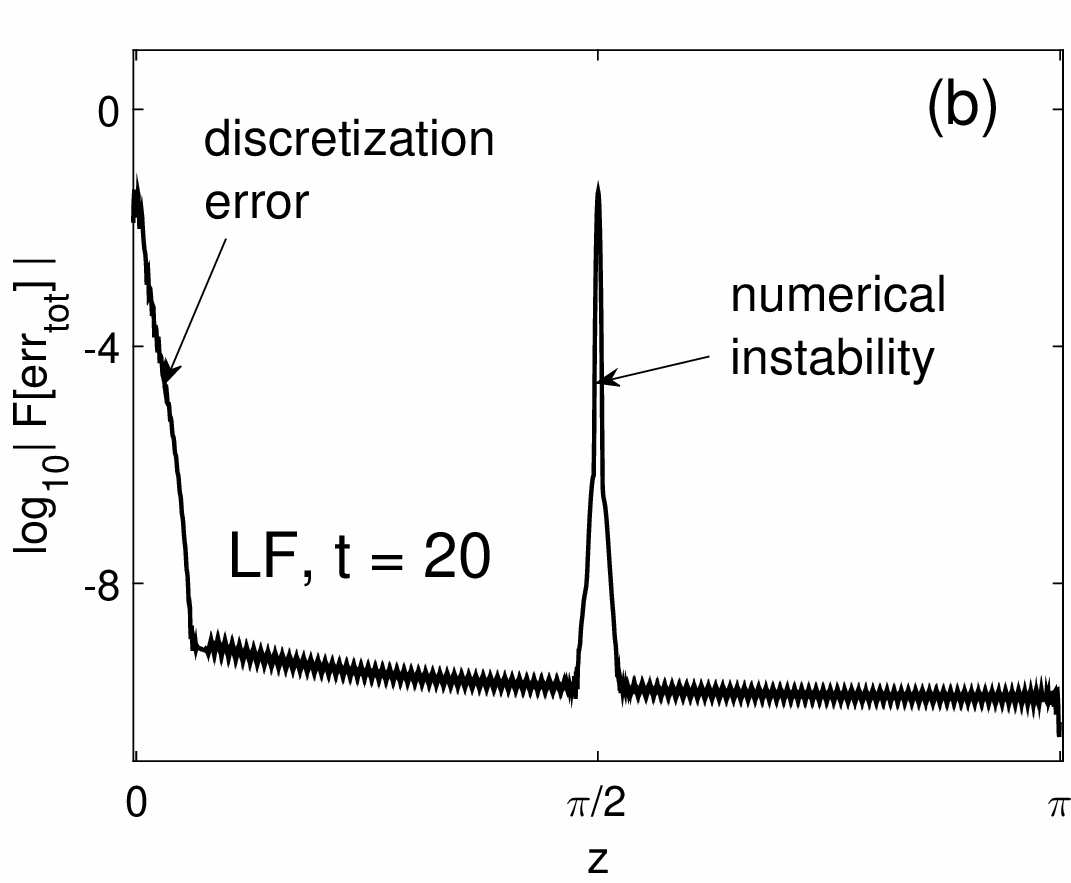}
\hspace{0.1cm}
\includegraphics[height=4.4cm,width=5.2cm,angle=0]{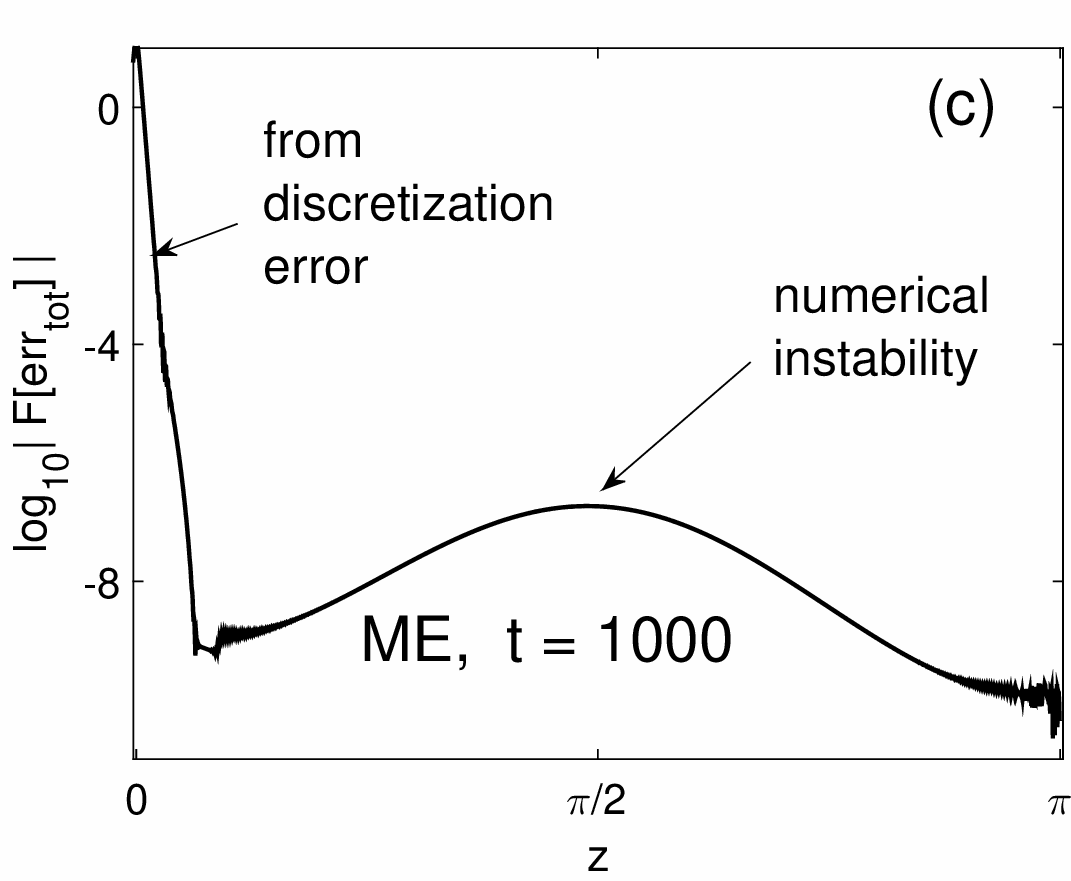}
\caption{
(a): Soliton solution \eqref{B03b} for $\Omega=0.7$. 
Note that Re$(V)$=Re$(U)$, Im$(V)$=$-$Im$(U)$. \ 
(b) and (c): Spectra of the total error \eqref{B04} obtained by the
MoC-LF at $t=20$ (b) \ and \ by the MoC-ME (c) at $t=1000$. 
See Section 7 and Appendix B for more details; in particular,
the origin of the peaks near $z=0$ in panels (b) and (c) is explained at 
the end of Appendix B. 
}
\label{fig_5}
\end{figure}

Let us now demonstrate the validity of qualitative conclusions
of our von Neumann analysis for a system whose linearization \eqref{e10a}
has a spatially dependent matrix ${\bf P}$. This system is the soliton
(see Fig.~\ref{fig_5}(a))
of the Gross--Neveu model \cite{GrossNeveu} 
in the relativistic field theory; its details are presented in Appendix B. 
It has received considerable attention in the past decade both from the analytical
and numerical perspectives: see \cite{2014_NLDE}, 
\cite{2012_BC}--\cite{2016_SciChi},
and references therein. 
The plots of the spectra of the numerical errors, obtained 
for this soliton with the MoC-LF
(Fig.~\ref{fig_5}(b)) and MoC-ME (Fig.~\ref{fig_5}(c)),
 are qualitatively similar to such plots in Sections 5 and 4, respectively.
In particular, the result in panel (b) agrees with our conclusion that the MoC-LF is 
unsuitable for simulations of a physically stable hyperbolic PDE system with
crossing characteristics:
the numerical instability at $k\approx \pi/(2h)$ will destroy the solution
soon after $t=20$. 
 In contrast, panel (c) shows that the MoC-ME would
be a suitable method for such systems, even for relatively long times.
Moreover, we also verified our conclusion (see \eqref{e26b} and 
Fig.~\ref{fig_2}(c)) 
that the growth rate of the most unstable Fourier harmonics,
with $k\approx \pi/(2h)$, of the MoC-ME scales as $O(h)$ (which is considerably
greater than the growth rate of that method's error for energy-preserving ODEs). 
To that end, we repeated the simulations for several values of $h=L/M$,
where $M=2^{11},\,2^{12},\,2^{13},\,2^{14}$ and computed $\gamma_{\rm meas}$
using \eqref{e23} and \eqref{B04}. The results fit the dependence
$\gamma_{\rm meas} = 6.4\cdot10^{-3}\cdot (2^{11}/M)\cdot\ln(10)$, which is
equivalent to $\gamma_{\rm meas}=O(h)$, to two significant figures.

\section*{Acknowledgement}

This work was supported in part by the NSF grant
 DMS-1217006.

\section*{Appendix A: \ A broader class of PDE systems with constant coefficients}

The class of models that we referred to in Section 7 generalizes Eqs.~\ref{e05}:
\be
\vSpm_{\;t} \pm \vSpm_{\;x} = \vSpm \times {\bf \hat{J}_c} \vSmp + 
                       \vSpm \times {\bf \hat{J}_s} \vSpm\,,
\label{A01}
\ee
where ${\bf \hat{J}_c}$ and ${\bf \hat{J}_s}$ are matrices accounting,
respectively, for cross- and self-interaction among components of the
Stokes vectors $\vSpm$. This class of models describes propagation of
electromagnetic waves in optical fibers with various types of birefringence.
The model considered in the main text is a special case of model \eqref{A01}
with 
\bsube
\be
{\bf \hat{J}_c} = \alpha \;{\rm diag}(1,-1,-2), \qquad
{\bf \hat{J}_s} = (2-3\alpha)\,{\rm diag}(0,0,1),
\label{A02a}
\ee
which corresponds to a highly spun birefringent fiber \cite{Wabnitz09}.
The parameter $\alpha$, which in model \eqref{e05} is set to $2/3$, 
accounts for the ellipticity of the fiber's core. Two other models are:
that of a randomly birefringent fiber \cite{Kozlov11}, with
\be
{\bf \hat{J}_c} = {\rm diag}(-1,1,-1), \qquad
{\bf \hat{J}_s} = \mathcal{O},
\label{A02b}
\ee
and that of an isotropic (i.e., non-birefringent) fiber \cite{Pitois01}:
\be
{\bf \hat{J}_c} = {\rm diag}(-2,-2,0), \qquad
{\bf \hat{J}_s} = {\rm diag}(-1,-1,0).
\label{A02c}
\ee
\label{A02}
\esube
Let us stress 
that only the above forms of matrices ${\bf \hat{J}_{c,\,s}}$ correspond
to physically meaningful situations in the context of birefringent fibers, 
and thus it would not
make sense to consider model \eqref{A01} with arbitrary ${\bf \hat{J}_{c,\,s}}$.

Models \eqref{A01}, \eqref{A02} each have many stationary constant solutions.
We will consider only six of them for each of those three models:
\be
S_j^+ = 1, \quad S_j^-=\pm S_j^+ \qquad 
\ba{l} \mbox{for one of $j=1,2$, or $3$, \ with} \\ 
  \mbox{the other two components of $\vSpm$ being $0$.} \ea 
\label{A03}
\ee
For brevity, we will refer to these solutions as $(j\pm)$, 
where $j$ and $\pm$ correspond to the particular choice of 
the component and the sign in \eqref{A03}.
For example, solution \eqref{e08} of the main text corresponds to
$(2-)$ in \eqref{A03}.

Of these 18 systems, the following 7 are stable on the infinite line
(in the sense specified before Eq.~\eqref{e12}):
model \eqref{A02a} with solutions $(1+)$, $(2-)$, $(3-)$; \ 
model \eqref{A02b} with solutions $(1-)$, $(2+)$, $(3-)$; \ 
model \eqref{A02c} with solution $(3+)$. \ 
The following 3 systems exhibit instability for $k=0$: \ 
model \eqref{A02a} with solutions $(1-)$, $(2+)$ and
model \eqref{A02c} with solution $(3-)$. 
(In other words, $\lambda_{\bf P}\in \mathbb{R}\setminus\{0\}$
for these systems.) The remaining 8 systems exhibit instability
for perturbations with $k\neq 0$.

\section*{Appendix B: \ The soliton solution of the Gross--Neveu model}

In physical variables, this model has the form \cite{GrossNeveu, 2014_NLDE}:
\be
\ba{l}
i( \psi_t + \chi_x ) + ( |\psi|^2 - |\chi|^2 )\psi - \psi = 0, \vspace{0.1cm}\\
i( \chi_t + \psi_x ) + ( |\chi|^2 - |\psi|^2 )\chi + \chi = 0 .
\ea
\label{B01}
\ee
%
A change of variables $u=(\psi+\chi)/\sqrt{2}$, $v=(\psi-\chi)/\sqrt{2}$
transforms \eqref{B01} to the form \eqref{e01}:
\be
\ba{l}
u_t + u_x = i(\,|v|^2 u + v^2u^*\,) - iv, \vspace{0.1cm} \\
v_t - v_x = i(\,|u|^2 v + u^2 v^*\,) -iu.
\ea
\label{B02}
\ee
The standing 
soliton solution of this system is (see, e.g., \cite{2014_NLDE} and references
therein):
\bsube
\be
\{u,v\} \,=\, \{U(x),V(x)\}\,\exp[-i\Omega t], 
\qquad \Omega \in (0,1);
\label{B03a}
\ee
\be
\{U(x),V(x)\}\, = \, \sqrt{1-\Omega}\,
        \frac{ \cosh(\beta x) \pm i\mu\,\sinh(\beta x) }
				{ \cosh^2(\beta x) - \mu^2\,\sinh^2(\beta x) };
\label{B03b}
\ee
\label{B03}
\esube
with $\beta=\sqrt{1-\Omega^2}$ and $\mu=\sqrt{(1-\Omega)/(1+\Omega)}$. 
The physical stability of this soliton for sufficiently small $\Omega$
has been an issue of recent controversy
\cite{2014_NLDE, 2015_NLDE_numerics} between analytical and numerical results.
However, for $\Omega$ sufficiently away from $0$, the soliton has been found to
be stable by both analysis and numerics (see references in \cite{2015_NLDE_numerics}). 
In Fig.~\ref{fig_5}(a) we show this solution for $\Omega=0.7$.

The linearized system \eqref{B02} on the background of soliton \eqref{B03} has
the form \eqref{e10a}, where the entries of matrix ${\bf P}$ are localized
functions of $x$. Therefore, a von Neumann analysis cannot be rigorously applied
to it. However, one can expect that its predictions could be valid qualitatively,
based on the principle of frozen coefficients. To demonstrate that this is indeed 
the case, we simulated system \eqref{B02} with the initial condition 
shown in Fig.~\ref{fig_5}(a), to which a small white noise was added. 
The length of the computational domain was $L=64$. The spectra
of the numerical
errors obtained by schemes \eqref{e27}, \eqref{e28} (MoC-LF) and 
\eqref{e24}, \eqref{e15} (MoC-ME) for $h=L/2^{12}\approx 0.016$
are shown in Figs.~\ref{fig_5}(b,c). The numerical error was defined similarly to
\eqref{e20}:
\be
{\rm err}_{\rm tot} = \left( \sum_{m=0}^M 
   \left| u_m^n - U(x_m)e^{-i\Omega t_n} \right|^2 + 
	 \left| v_m^n - V(x_m)e^{-i\Omega t_n} \right|^2 \, \right)^{1/2}.
\label{B04}
\ee
Note that the peaks near $z=0$, seen in Figs.~\ref{fig_5}(b,c), 
do {\em not} correspond to any numerical instability.
In panel (b) that peak is caused directly by the discretization error of the scheme.
In panel (c), a much higher such peak is also caused by the discretization error,
but indirectly. Indeed, a slight 
error in the propagation constant $\Omega$, which inevitably occurs due to
a limited accuracy of the numerical scheme, causes $\{u_m^n,\,v_m^n\}$ in \eqref{B04}
to differ from their respective $\{U(x_m),\,V(x_m)\}\,\exp[-i\Omega t_n]$ by 
$O(1)$ for $t_n\gg 1$.



\begin{thebibliography}{99}

%
\bibitem{GKO}
B.~Gustafsson, H.-O.~Kreiss, J.~Oliger, \ 
Time-dependent problems and difference methods, \ 2nd Ed., 
John Wiley \& Sons, Hoboken, NJ, 2013; \ Sec.~7.4.
%
\bibitem{stab_3DMoC}
Yu.Ya.~Mikhailov, \
Stability of some numerical schemes of the three-dimensional method of characteristics, \
USSR Comp.~Math.~Math.~Phys. 8 (1968) 312--315.
%
\bibitem{Australian08}
A.C.~Zecchin, A.R.~Simpson, M.F.~Lambert, \ 
von Neumann stability analysis of a method of characteristics visco-elastic
pipeline model, \ 10th Int.~Conf.~on Pressure Surges (Edinburgh, Scotland),
2008. 
%
\bibitem{ZM87}
V.E.~Zakharov, A.V.~Mikhailov, \ 
Polarization domains in nonlinear media, \ 
JETP Lett. 45 (1987) 349--352.
%
\bibitem{Pitois01}
S.~Pitois, G.~Millot, S.~Wabnitz, \ 
Nonlinear polarization dynamics of counterpropagating waves in an isotropic
optical fiber: theory and experiments, \
J.~Opt.~Soc.~B 18 (2001) 432--443.
%
\bibitem{Wabnitz09}
S.~Wabnitz, \ 
Chiral polarization solitons in elliptically birefringent spun optical fibers, \
Opt.~Lett. 34 (2009) 908--910.
%
\bibitem{Kozlov11}
V.V.~Kozlov, J.~Nuno, S.~Wabnitz, \ 
Theory of lossless polarization attraction in telecommunication fibers, \ 
J.~Opt.~Soc.~B 28 (2011) 100--108.
%
\bibitem{Zthesis}
Z.~Deng, \ 
Uncommon numerical instability in the method of characteristics
applied to hyperbolic equations, \ 
M.S.~Thesis, University of Vermont, 2016. 
%
\bibitem{Boyd}
R.W.~Boyd, \ Nonlinear optics, \ Academic, San Diego, 1992.
%
\bibitem{SBS_Kaup}
D.J.~Kaup, \ 
The first-order perturbed SBS equations, \ 
J.~Nonlin.~Sci. 3 (1993) 427--443.
%
\bibitem{SBS_JQE}
C.E.~Mungan, S.D.~Rogers, N.~Satyan, J.O.~White, \ 
Time-dependent modeling of Brillouin scattering in optical fibers excited
by a chirped diode laser, \ 
IEEE J.~Quant.~Electron. 48 (2012) 1542--1546. 
%
\bibitem{SRS_1975}
F.Y.F.~Chu, A.C.~Scott, \ 
Inverse scattering transform for wave-wave scattering, \ 
Phys.~Rev.~A 12 (1975) 2060--2064. 
%
\bibitem{RKB_1979}
D.J.~Kaup, A.~Rieman, A.~Bers, \ 
Space-time evolution of nonlinear three-wave interactions. I. Interactions
in an homogeneous medium, \ 
Rev.~Mod.~Phys. 51 (1979) 275--310. 
%
\bibitem{SGH_Kaup}
D.J.~Kaup, \ 
Simple harmonic generation: an exact method of solution, \ 
Stud.~Appl.~Math. 59 (1978) 25--35.
%
\bibitem{Ibragimov}
E.~Ibragimov, A.~Struthers, D.J.~Kaup, \ 
Parametric amplification of chirped pulses in the presence of a large
phase mismatch, \ 
J.~Opt.~Soc.~Am.~B 18 (2001) 1872--1876. 
%
\bibitem{CJM_quant1}
C.J.~McKinstrie, L.~Mejling, M.G.~Raymer, K.~Rottwitt, \ 
Quantum-state-preserving optical frequency conversion and 
pulse reshaping by four-wave mixing, \ 
Phys.~Rev.~A 85 (2012) 053829. 
%
\bibitem{CJM_quant2}
D.V.~Reddy, M.G.~Raymer, C.J.~McKinstrie, \ 
Efficient sorting of quantum-optical wave packets by 
temporal-mode interferometry, \ 
Opt.~Lett. 39 (2014) 2924--2927.
%
\bibitem{CJM_photon}
C.J.~McKinstrie, D.S.~Cargill, \ 
Simultaneous frequency conversion, regeneration and 
reshaping of optical signals, \ 
Opt.~Expr. 20 (2012) 6881--6886.
%
\bibitem{CJM_plasma1}
C.J.~McKinstrie, E.J.~Turano, \ 
Spatiotemporal evolution of parametric instabilities driven by short
laser pulses: One-dimensional analysis, \ 
Phys.~Plasmas 3 (1996) 4683--4696. 
%
\bibitem{CJM_plasma2}
C.J.~McKinstrie, V.A.~Smalyuk,  R.E.~Giacone, H.X.~Vu, \ 
Power exchange between crossed laser beams and the associated frequency cascade, \ 
Phys.~Rev.~E 55 (1997) 2044--2047. 
%
\bibitem{FBGbook}
R.~Kashyap, \ 
Fiber Bragg Gratings, \ Academic, Burlington, MA, 2010; Chap.~4.
%
\bibitem{2002_CW}
J.E.~Sipe, C.M.~de Sterke, B.J.~Eggleton, \ 
Rigorous derivation of coupled mode equations for short,
high-intensity grating-coupled, co-propagating pulses, \ 
J.~Mod.~Opt. 49 (2002) 1437–-1452.
%
\bibitem{2014_ChowdAtai}
S.A.M.S.~Chowdhury, J.~Atai, \ 
Stability of Bragg grating solitons in a semilinear dual core system
with dispersive reflectivity, \ 
IEEE J.~Quant.~Electron. 50 (2014) 458--465. 
%
\bibitem{Thirring}
W.E.~Thirring, \ 
A soluble relativistic field model, \ 
Ann.~Phys. 3 (1958) 91--112.
%
\bibitem{GrossNeveu}
D.J.~Gross, A.~Neveu, \ 
Dynamical symmetry breaking in asymptotically free field theories, 
Phys.~Rev.~D 10 (1974) 3235--3253.
%
\bibitem{2006_ChugunovaPeli}
M.~Chugunova, D.~Pelinovsky, \ 
Block-diagonalization of the symmetric first-order coupled-mode system, \ 
SIAM J.~Appl.~Dyn.~Syst. 5 (2006) 55–-83.
%
\bibitem{2014_NLDE}
S.~Shao, N.R.~Quintero, F.G.~Mertens, F.~Cooper, A.~Khare, A.~Saxena, \ 
Stability of solitary waves in the nonlinear Dirac equation with arbitrary
nonlinearity, \ 
Phys.~Rev.~E 90 (2014) 032915.
%
\bibitem{AcevesWabnitz}
A.B.~Aceves, S.~Wabnitz, \ 
Self-induced transparency solitons in nonlinear refractive periodic media, \ 
Phys.~Lett.~A 141 (1989) 37--42. 
%
\bibitem{1992_RTW}
M.~Romagnoli, S.~Trillo, S.~Wabnitz, \ 
Soliton switching in nonlinear couplers, \ 
Opt.~Quantum Electron. 24 (1992) S1237--1267. 
%
\bibitem{Assemat12}
E.~Assemat, A.~Picozzi, H.-R.~Jauslin, D.~Sugny, \ 
Hamiltonian tools for the analysis of optical polarization control, \ 
J.~Opt.~Soc.~B 29 (2012) 559--571.
%
\bibitem{Griffiths_book}
D.F.~Griffiths, D.J.~Higham, \ 
Numerical methods for ordinary differential equations, \
Springer-Verlag, London, 2010; Chaps. 13 and 15.
%
\bibitem{2012_BC}
G.~Berkolaiko, A.~Comech, \ 
On spectral stability of solitary waves of nonlinear Dirac equation in 1D, \ 
Math.~Model.~Nat.~Phenom. 7 (2012) 13--31.
%
\bibitem{2012_PS}
D.E.~Pelinovsky, A.~Stefanov, \ 
Asymptotic stability of small gap solitons in the nonlinear Dirac equations, \ 
J.~Math.~Phys. 53 (2012) 073705.
%
\bibitem{2016_BC}
N.~Boussa\"id, A.~Comech, \ 
On spectral stability of the nonlinear Dirac equation, \ 
J.~Funct.~Anal. 271 (2016) 1462--1524.
%
\bibitem{2010_CPC}
F.~de la Hoz, F.~Vadillo, \ 
An integrating factor for nonlinear Dirac equations, \ 
Comp.~Phys.~Commun. 181 (2010) 1195--1203. 
%
\bibitem{2013_JCP}
J.~Xu, S.~Shao, H.~Tang, \ 
Numerical methods for nonlinear Dirac equation, \ 
J.~Comp.~Phys. 245 (2013) 131--149.
%
\bibitem{2015_NLDE_numerics}
J.~Cuevas-Maraver, P.G.~Kevrekidis, A.~Saxena, F.~Cooper, F.G.~Mertens, \ 
Solitary waves in the nonlinear Dirac equation at the continuum limit: 
Stability and dynamics, \ 
In: Ord.~Part.~Diff.~Eqs., Nova Science, Boca Raton, 2015; Chap.~4.
%
\bibitem{2016_SciChi}
W.Z.~Bao, Y.Y.~Cai, X.W.~Jia, J.~Yin, \ 
Error estimates of numerical methods for the nonlinear Dirac equation 
in the nonrelativistic limit regime, \ 
Science China Math. 59 (2016) 1461--1494.








\end{thebibliography}
\end{document}